\documentclass[11pt]{amsart}
\usepackage{amsmath}
\usepackage{soul}
\usepackage{url}
\usepackage{hyperref}

\usepackage[margin=3.5cm]{geometry}

\makeatletter
\newcommand{\addresseshere}{%
  \enddoc@text\let\enddoc@text\relax
}
\makeatother

\usepackage{algorithm,algpseudocode,float}

\usepackage{amssymb}
\usepackage{enumerate}
\usepackage{amsmath}
\usepackage{amsthm}
\usepackage{tikz} 
\usepackage{float}
\usepackage{mathtools}
\usepackage{graphicx}
\usepackage{pgfplots}
\pgfplotsset{compat=1.18}
\allowdisplaybreaks

\newcommand{\pn}{\mathbb{P}_n}
\newcommand{\bpn}{\mathbb{BP}_n}
\newcommand{\pmn}{\mathbb{P}_m(n)}
\newcommand{\flop}{{r}^{-1}}
\newcommand{\AlCyc}{\textsc{AlCyc}}
\newcommand{\AlGra}{\textsc{AlGra}}

\newcommand*\circled[1]{%
  \tikz[baseline=(char.base)]\node[anchor=south west, draw,rectangle, rounded corners, inner sep=2pt, minimum size=7mm,
    text height=2mm](char){\ensuremath{#1}} ;}

\usepackage{pgf,tikz}\usepackage{mathrsfs}\usetikzlibrary{arrows}

\usetikzlibrary{matrix,fit, arrows, shapes}

\numberwithin{equation}{section}
\newtheorem{thm}{Theorem}[section]
\newtheorem{prop}[thm]{Proposition}
\newtheorem{cor}[thm]{Corollary}

\newtheorem{rem}[thm]{Remark}
\newtheorem{conj}[thm]{Conjecture}

\theoremstyle{definition}
\newtheorem{defn}{Definition}[section]
\newtheorem{ex}[defn]{Example}

\begin{document}

\author{Sa\'ul A. Blanco and Charles Buehrle}%\corref{cor}}
\title[Bounds on the genus for prefix-reversal graphs]{Bounds on the genus for 2-cell embeddings of prefix-reversal graphs}

\address{Department of Computer Science, Indiana University, Bloomington, IN 47408}
\email{sblancor@indiana.edu}
\address{ Department of Mathematics, Physics, and Computer Studies, Notre Dame of Maryland University, Baltimore, MD 21210}
\email{cbuehrle@ndm.edu}

\date{September 23, 2024}

\begin{abstract}
In this paper, we provide bounds for the genus of the pancake graph $\mathbb{P}_n$, burnt pancake graph $\mathbb{BP}_n$, and undirected generalized pancake graph $\mathbb{P}_m(n)$. Our upper bound for $\mathbb{P}_n$ is sharper than the previously-known bound, and the other bounds presented are the first of their kind. Our proofs are constructive and rely on finding an appropriate rotation system (also referred to in the literature as Edmonds' permutation technique) where certain cycles in the graphs we consider become boundaries of regions of a 2-cell embedding. A key ingredient in the proof of our bounds for the genus $\mathbb{P}_n$ and $\mathbb{BP}_n$ is a labeling algorithm of their vertices that allows us to implement rotation systems to bound the number of regions of a 2-cell embedding of said graphs. All of our bounds are asymptotically tight; in particular, the genus of $\mathbb{P}_m(n)$ is $\Theta(m^nnn!)$ for all $m\geq1$ and $n\geq2$.

\end{abstract}

\bibliographystyle{plain}

\maketitle

\section{Introduction}

In graph theory, determining if a graph is \emph{planar}, that is, if one can draw the graph in the plane (or on a sphere) preventing its edges from intersecting, has been of interest since Euler's time. Indeed, it is a classical result, referred to as the \emph{Euler polyhedral formula}, that if $G$ is a planar graph, then
\[v-e+r=2,\]
where $v$ is the number of vertices of $G$, $e$ is the number of edges of $G$, and $r$ is the number of regions into which $G$ divides the plane. In the 1930s, Kuratowski~\cite{K30} and Wagner~\cite{W37} provided a complete characterization of what graphs are planar in terms of the complete graph $K_5$ and the complete bipartite graph $K_{3,3}$.

By a similar token, once a graph is known to not to be planar, the next natural question is determining what is the orientable surface of minimum genus (in the topological sense) in which $G$ can be embedded (drawn without intersections). Indeed, the terms \emph{toroidal graph, double-toroidal graph} and \emph{pretzel graph} are used to refer to graphs where the minimum-genus surface they can be embedded on are the torus (genus 1), double torus (genus 2), or a triple torus (genus 3), respectively, for reference see West~\cite[Page 266]{West00}.

As a computational problem, determining if a given a graph $G$ has genus $g$ or less, where $g$ is a positive integer, is an NP-complete problem by Thomassen~\cite{T89}. Moreover, the genus of a graph is tied to the computational complexity of certain classical graph problems. In particular, the fixed-parameter tractability, subexponential time computability, and approximability of \textsc{Independent Set}, \textsc{Vertex Cover}, and \textsc{Dominating Set} depends on the genus of the corresponding graph, see Chen, et al.~\cite{C07}.

Another interesting appearance of the genus is in \emph{proof-labeling schemes}, first introduced in Korman, Kutten, and Peleg~\cite{KKP05} whose formal definition is beyond the scope of this paper. Intuitively, a proof-labeling scheme is a protocol that certifies locally if a certain property holds on a network. By ``locally," one means that each vertex $v$ in the network has a partial view of the entire network, such as the set of vertices adjacent to $v$. If a graph satisfies a property $\mathcal{P}$, then every vertex has to certify $\mathcal{P}$, whereas if the graph does not satisfy the property then at least one vertex fails to certify $\mathcal{P}$. Planar graphs and graphs of bounded genus have relatively straight forward proof-labeling schemes (see Feuilloley, et al.~\cite{Fetal20, Fetal23}). If a graph $G$ can be embedded on a surface of genus $g$, then $G$ has a proof-labeling whose complexity depends on $g$ (see Esperet and Lévêque as well as Feuilloley, et al.~\cite{EL22,Fetal23}).

A family of graphs that is ubiquitous in computer science and discrete mathematics is the family of \emph{pancake graphs} or \emph{prefix-reversal graphs}. These graphs have applications to parallel computing and bioinformatics (see, for example, Hannenhalli and Pevzner, Hayes, et al., Kanevsky and Feng, and Lakshmivarahan, Jwo, and Dhall~\cite{HannenPev,Haynes2008,KF95,LJD93}). These graphs arise from the classical \emph{pancake problem} of Dweighter~\cite{Dweighter75} where one wishes to sort a stack of pancakes of different diameters by flipping some of the pancakes at a time utilizing a chef's spatula, and repeating this process. Several properties of these graphs are understood, such as their being regular, vertex-transitive, and having a recursive structure. Regarding the cycles in these graphs, one can find cycles of all length from the \emph{girth} (length of the shortest cycle) to a Hamiltonian cycle, see the authors' other work (with Patidar) and Kanevsky and Feng~\cite{BBP19, BB23, KF95}. Regarding the genus of the pancake graph, Nguyen and Bettayeb~\cite{NB09} provided the first known lower and upper bounds. In this paper, we improve their upper bound, and  provide lower and upper bounds for other related pancake graphs. 

\subsection{Our contribution} Our main contributions are the following results. 
\begin{enumerate}[(i)]
    \item We provide a tighter upper bound to that in Nguyen and Bettayeb~\cite{NB09} for the pancake graph $\pn$. We show in Theorem~\ref{t:pnupperbound} that if $n>3$, then
    \[
    \gamma(\pn) \leq n!\left(\frac{3n-10}{12}\right) + 1.
    \]
    \item We provide lower and upper bounds for the genus of the burnt pancake graph $\bpn$. We show in Corollary~\ref{cor:bpnlowerbound} and Theorem~\ref{t:bpnupperbound} that if $n>2$, then
    \[
     2^{n - 4} (3 n - 8) n! + 1\leq \gamma(\bpn) \leq 2^{n-4}(4n-9)n! + 1.
    \]
    \item We provide lower and upper bounds for the genus of the undirected generalized pancake graph $\pmn$. Indeed, in Corollary~\ref{cor:pmnlowerbound} and Theorem~\ref{t:pmnupperbound} we establish that if $n>1$, then 
    \[
        \gamma(\pmn)\geq \begin{cases}
        \frac{1}{2}m^{n - 1} ((m-2) n - m) n! + 1, & m \in \{3,4,5\}\\
        \frac{1}{6}m^n (2n-3) n! + 1, & m \geq 6,
    \end{cases}
    \] and that
        if $m \geq 3$ and $n \geq 2$, then
    \[
        \gamma(\pmn)\leq
\begin{cases}
    \frac{1}{2}m^{n-1} (mn-m-n)n!+1,&\text{ if }m\text{ is even}\\
    \frac{1}{2}m^{n-1} (2 m n - 2 m - n-1) n!+1,&\text{ if }m\text{ is odd.}
\end{cases} 
    \] 
    \item As a consequence of (i)-(iii) above, $\gamma(\pmn)$ is $\Theta(m^nnn!)$\footnote{$f(n)$ is said to be $\Theta(g(n))$ if there exist positive real numbers $C_1,C_2$ and positive integer $n_0$ such that $C_1|g(n)|\leq |f(n)|\leq C_2|g(n)|$ for all $n>n_0$.} for all $m\geq1$ and $n\geq2$.
    \item In the particular case of  $\mathbb{P}_m(2)$, we are able to provide tighter bounds than the ones given in Theorem~\ref{t:bpnupperbound}. In Remark~\ref{r:pm2} we show that
     \[
        \gamma(\mathbb{P}_m(2)) \leq \begin{cases}
            m^2-3m+1, &\text{ if } m \text{ odd and } m=3j\\
            m^2-2m+1, &\text{ if } m \text{ odd and } m\neq3j\\
            m^2-(7/2) m+1, &\text{ if } m \text{ even and } m=3j\\
            m^2-(5/2) m+1, &\text{ if } m \text{ even and } m\neq3j,\\
        \end{cases}
    \] and we furthermore establish that $\gamma(\mathbb{P}_3(2))=1$, or that $\mathbb{P}_3(2)$ is toroidal.
\end{enumerate}

\begin{rem}
    We remark that all our bounds are asymptotically tight as well. 
\end{rem}

\subsection{Organization of the paper} The remainder of this paper is organized as follows. In Section~\ref{sec:preliminaries} we present the necessary basic definitions and notation used throughout. These include formal definitions of the family of pancake graphs treated in this paper as well as rotation systems and 2-cell embedding of graphs. In particular, we explain how to obtain a 2-cell embedding of graphs using rotation systems (also known as Edmonds' permutation technique). In Section~\ref{sec:labelingpnbpn}, we give an algorithm, which we call $\AlGra$, standing for ``\textbf{al}ternating \textbf{gra}ph" labeling to label the vertices of the classical pancake and burnt pancake graphs in a way that certain cycles become regions in a 2-cell embedding given by an appropriate rotation system. In Section~\ref{sec:bounds}, we use the labeling given by $\AlGra$ and define rotation systems that yield a lower bound on the number of regions and thus give upper bounds for the genus of the pancake graph and burnt pancake graphs. We provide bounds for the genus of the undirected generalized pancake graph as well by defining another rotation system that takes advantage of the fact that in the generalized pancake graphs, prefix reversals are \emph{not} their own inverses. In addition, we also provide lower bounds for the genus of the burnt pancake and undirected generalized pancake graphs that depends on their girth. 

\section{Preliminaries}\label{sec:preliminaries}

\subsection{Graph embeddings}

We follow traditional graph theory convention and denote an undirected graph by $G=(V,E)$ where $V=V(G)$ is the set of vertices of $G$ and $E=E(G)\subset V\times V$ is the set of edges of $G$. An \emph{embedding} (``imbedding" is also an acceptable spelling) of $G$ is a realization of the graph where every $e=\{u,v\}\in E$ is represented by curve segments (that is, curves homeomorphic to $[0,1]$) whose endpoints are $u,v$, and where no such curve segments intersect except at the endpoints. The \emph{genus} of $G$, denoted by $\gamma(G)$, is defined as
\[
\min_S\{\gamma(S):\text{$G$ can be embedded in S}\}, 
\] where $S$ is an orientable surface and $\gamma(S)$ denotes the genus of $S$ in the topological sense. If we consider the complement of the curve segments that form the embedding of $G$ on a surface $S$, one obtains the \emph{regions} of the embedding. If each of the regions is homeomorphic to an open disk, we say that the embedding is a \emph{2-cell embedding} or a \emph{cellular embedding}. The boundary of a region corresponds to a \emph{cycle} or \emph{boundary cycle}, which is a sequence of edges where only the first and last vertex are equal. The Euler-Poincar\'e formula relates the number of vertices, edges, and regions with the genus of $G$. Indeed, if $G=(V,E)$ is a graph with a 2-cell embedding on an orientable surface, then it follows that
\[
|V|-|E|+r=2-2\gamma(G),
\] where $r$ is the number of regions of the embedding.

\subsection{Generalized symmetric group and prefix reversals}

Let $m,n$ be two positive integers. Let $C_m$ denote the \emph{cyclic group} of order $m$ and $S_n$ denote the \emph{symmetric group} of order $n!$. Consider the group $S(m,n):=C_m\wr S_n$, the wreath product of $C_m$ and $S_n$. $S(m,n)$ is usually referred to as the \emph{generalized symmetric group.} We shall denote the elements of $S(m,n)$ as a string of $n$ distinct characters taken from the set $[n]:=\{1,2,\ldots, n\}$, each of which has a superscript representing a ``sign" taken from the set $\{0\}\cup[m-1]$. For example, $3^04^21^22^35^1\in C_4\wr S_5$ has five characters, the elements of $[5]$, and four signs, the elements of $\{0\}\cup[3]$.

We define a pair of sets of \emph{generalized prefix reversals}. First, we define a \emph{flip} $r_i$, with $1\leq i\leq n$, to be the function $r_i:S(m,n)\to S(m,n)$ given as follows. If $\pi=\pi_1^{a_1}\pi_2^{a_2}\cdots \pi_n^{a_n}\in S(m,n)$, the flip reversal $r_i$ will reverse the first $i$ characters of $\pi$ and add one to each of the corresponding signs modulo $m$. In other words,
\[
r_i(\pi_1^{a_1}\pi_2^{a_2}\cdots \pi_n^{a_n})= \pi_i^{b_i}\pi_{i-1}^{b_{i-1}}\cdots 
\pi_1^{b_1}\pi_{i+1}^{a_{i+1}}\cdots \pi_n^{a_n}, 
\] where $b_j=a_j+1\mod m$ for $1\leq j\leq i$.
Second, we define a \emph{flop} $\flop_i$, with $1\leq i \leq n$, to be the function $\flop_i:S(m,n)\to S(m,n)$ given as follows. If $\pi=\pi_1^{a_1}\pi_2^{a_2}\cdots \pi_n^{a_n}\in S(m,n)$, then
\[
\flop_i(\pi_1^{a_1}\pi_2^{a_2}\cdots \pi_n^{a_n})= \pi_i^{c_i}\pi_{i-1}^{c_{i-1}}\cdots 
\pi_1^{c_1}\pi_{i+1}^{a_{i+1}}\cdots \pi_n^{a_n}, 
\] where $c_j=a_j-1\mod m$ for $1\leq j\leq i$.

\subsection{The cases \texorpdfstring{$S(1,n)$}{S(1,n)} and \texorpdfstring{$S(2,n)$}{S(2,n)}} 

The group $S(1,n)$ is isomorphic to $S_n$, and in this case the sign is not indicated. So $4^03^05^01^02^0$ is simply written as $43512$. In $S_n$, the prefix reversal $r_i$ (equivalently $\flop_i$) will simply reverse the first $i$ characters of a permutation and it is therefore an involution (i.e., each prefix reversal is its own inverse). Additionally, the action of $r_1$ in $S_n$ is trivial and therefore we omit it from the set of prefix reversals of $S_n$.

The group $S(2,n)$ is isomorphic to the \emph{hyperoctahedral group}, or the group of \emph{signed permutations}. In this context, the signs $\{0,1\}$ are simply written as $\{+,-\}$. When listing elements of $S(2,n)$ the ``$+$" is omitted and the ``$-$" is usually denoted by an underline. For example, $2^13^04^05^11^0$ is written as $\underline{2}34\underline{5}1$. The usual definition of the hyperoctahedral group (for reference see Bj\"{o}ner and Brenti~\cite[Chapter 8]{BjornerBrenti}) is the group of permutations $w$ on the set $[\pm n]:=\{-n,-(n-1),\ldots, -2,-1,1,2,\ldots, n\}$ with $w(-i)=-w(i)$ for $1\leq i\leq n$. Moreover, the hyperoctahedral group on $[\pm n]$ is denoted by $B_n$. Notice that in the case of $B_n$, a prefix reversal $r_i$ will reverse the first $i$ characters of a signed permutation, and will flip each of the signs of the first $i$ characters from ``$+$" to ``$-$" and vice versa. Therefore, the prefix reversals in $B_n$ are involutions. 

The following remark is worth making since the cases $S(1,m)$ and $S(2,n)$ are of particular interest in the literature. 

\begin{rem}
    $S(1,n)$ is isomorphic to $S_n$ and $S(2,n)$ is isomorphic to $B_n$. Moreover, there are $n$ prefix reversals for $B_n$, labeled $r_i$ with $1\leq i\leq n$ and there are $n-1$ prefix reversals for $S_n$, labeled $r_i$ with $2\leq i\leq n$. Additionally the prefix reversals for $S(1,n)\cong S_n$ and $S(2,n)\cong B_n$ are all involutions, that is, $r_i=r^{-1}_i$ for $1\leq i\leq n$. 
\end{rem}

\subsection{Generalized pancake graphs \texorpdfstring{$P(m,n)$}{P(m,n)}, related graphs, and their recursive structure}

The family of \emph{generalized pancake graphs} (also referred to as \emph{$k$-sided pancake network} in Cameron, Sawada, and Williams~\cite{CSW21}) consists of the Cayley graphs of $S(m,n)$ using the prefix reversals of $S(m,n)$ as generators utilizing the composition operation. In other words, the vertex set of a generalized pancake graph of $S(m,n)$, denoted by $P(m,n)$, is the graph with vertex set $S(m,n)$ and if $v_1,v_2\in S(m,n)$ then $(v_1,v_2)$ is an edge of $P(m,n)$ if and only if $v_2=r_i(v_1)$ for some $1\leq i\leq n$. There are three families of graphs that are of special interest, and whose genus are the focus of this paper. 
\begin{description}
    \item[Undirected generalized pancake graph, $\pmn$] This graph arises by ignoring the directions of the edges in $P(m,n)$. The cycle structure of $\pmn$ has been studied by the authors in~\cite{BB23} where they provided the girth (i.e., the length of the shortest cycle) of $\pmn$ and proved that these graphs are weakly pancyclic (i.e., $\pmn$ has all cycles of length from its girth to a Hamiltonian cycle.) 
    \item[The pancake graph, $\pn$] The special case $P(1,n)$. The name of $\pn$ comes from the classical \emph{pancake problem}, first proposed in Dweighter~\cite{Dweighter77}, which is equivalent to finding the diameter of $\pn$. 
    \item[The burnt pancake graph, $\bpn$] The special case $P(2,n)$. The name of $\bpn$ comes from the \emph{burnt pancake problem}, first proposed by Gates and Papadimitriou in~\cite{GatesPapa}.
\end{description}

We label the edges of our graphs by the prefix reversal used to get from one vertex to another. 

The generalized pancake graphs $P(m,n)$ and the special cases mentioned above have a recursive structure. To help us describe said structure, we are going to use $\cdot$ to denote concatenation of two strings. For example, $12\cdot 453$ produces $12453\in S_5$, $\underline{1}3\underline{5}\cdot 2\underline{4}$ produces $\underline{1}3\underline{5} 2\underline{4}\in B_5$, and $1^22^04^1 \cdot 5^23^0$ produces $1^22^04^15^23^0 \in S(3,5)$. Furthermore, we use $\ell(\cdot)$ to denote the \emph{length function} on strings. By abuse of notation, we use $e=1^02^0 \cdots n^0$ to denote the identity element of $S(m,n)$ since the values of $m,n$ will be clear from context.  

The recursive structure of our graphs is exhibited readily by looking at the induced subgraphs by the set of vertices sharing the same character and sign. For example, if the reader takes a look at Figure~\ref{f:P4}, we can see four subgraphs that are isomorphic to $\mathbb{P}_3$ (all isomorphic to a 6-cycle) by considering the sets

\[
A_i=\{w\in \mathbb{P}_4: w=w'\cdot i\}\text {, with } i\in[4].
\]

In general, if one takes $w\in S(m,n)$ and consider a \emph{suffix} $s$ of $w$, that is, if we write $w=w'\cdot s$, then the subgraph of $P(m,n)$ induced by the set
\[
A_s=\{w\in S(m,n): w=w'\cdot s\}
\]
is isomorphic to $P(m,k)$, with $k=n-\ell(s)$. We denote these subgraphs by $P(m,k)\cdot s$.  

In practice, we will be dealing mostly with $\pn,\bpn$, and $\pmn$ with $m>2$. These graphs inherit an analogous recursive structure from $P(m,n)$. For example, $\mathbb{P}_2\cdot 34$ corresponds to the edge $(1234,2134)\in E(\mathbb{P}_4)$.

\subsection{Convention relating cycles}\label{sec:prelim-cycles} By a \emph{cycle} in a graph $G$ of length $\ell$, we mean a sequence of vertices and edges of the form $v_1e_1\cdots v_{\ell}e_{\ell}v_{\ell+1}$ with $e_i=\{v_i,v_{i+1}\}$, $1\leq i\leq \ell$, where all the edges and vertices are different except for $v_1=v_{\ell+1}$. Moreover, if $C=v_1e_1\cdots v_{\ell}e_{\ell}v_{\ell+1}$ is a cycle, we sometimes refer to it as an \emph{$\ell$-cycle} when we wish to emphasize the length of $C$. Since the graphs we deal with in this paper are all simple, meaning that every edge has two different endpoints and that any pair of distinct vertices are the endpoints of at most one edge, we will omit the vertices from the notation so $C=v_1e_1\cdots v_{\ell}e_{\ell}v_{\ell+1}$ becomes $C=e_1\cdots e_{\ell}$.

Furthermore, if we can traverse a cycle $C$ by alternating the labels of two edges, say $e_1,e_2$, then we say that $C$ is a $\{e_1,e_2\}$-cycle.

\subsection{Rotation systems/Edmonds' permutation technique}

Given a graph $G$, a rotation $p_v$ of $v\in G$ is a cyclic permutation of the edges incident with $v$. A \emph{rotation system} $p$ is an assignment of a rotation to each vertex of $G$. %Formally, if $G$ has $n$ vertices $v_1,v_2,\ldots,v_n$, and $N(v_i)$ denote the sets of neighbors of $v_i$, $1\leq i\leq n$, then a rotation system is $p=\{p_{v_1},p_{v_2},\ldots, p_{v_n}\}$, where each $p_{v_i}$ is a cyclic permutation of length $|N(v_i)|$ of the elements of $N(v_i)$.
In our proofs we will be utilizing the following result from Wan~\cite{W22}.

\begin{prop}(\cite[Proposition 1]{W22})\label{prop:one}
    Given a simple, undirected graph $G$, there exists a bijection between the rotation systems and the orientable 2-cell embeddings of $G$. 
\end{prop}

The ideas for this result can be traced back to the 1800s, for examples see Dyck and Heffter~\cite{D88,H91,H98} and were used in 1965 in Ringel~\cite{R65}. Since then, it was independently discovered by Edmonds in his M.S. thesis~\cite{E60} and a more algorithmic treatment is provided in Youngs~\cite{Y63}.  Due to Edmonds' independent discovery, Proposition~\ref{prop:one}, or more precisely the algorithmic version provided in Youngs~\cite{Y63}, is also known as \emph{Edmonds' permutation technique}. We will employ Edmonds' permutation technique to derive the main results of this paper. A similar approach was used in White~\cite{W72} to determine the genus of certain Cayley graphs.

We will need some notation. Let $G=(V,E)$ be a simple graph and for $x\in V$, we let $N(x)$ be the set of \emph{neighbors} of $x$, namely, $N(x)=\{y\in V: \{x,y\}\in E\}$. For $x\in V$, define $p_x:N(x)\to N(x)$ to be a cyclic permutation of length $|N(x)|$. Let $E^o(G)=\{(x,y),(y,x):\{x,y\}\in E\}$; in other words, $E^o$ is the set of all ordered pairs (the ``$o$" superscript stands for ``ordered") that can be formed with edges from $G$, and $p:E^o\to E^o$ given by $p(x,y)=(y,p_y(x))$. Proposition~\ref{prop:one} states that any 2-cell embedding of $G$ uniquely determines the set of cyclic permutations $\{p_x\}_{x\in V}$, and that furthermore given the set $\{p_x\}_{x\in V}$, the orbits of $p$ determine a 2-cell embedding of $G$. Intuitively, if a surface is oriented, there is a unique cyclic ordering of each of the vertices' neighbors representing the ``clockwise" direction. 

We denote the orbit of an edge $(x,y)\in E^o$ by 
\[\{p^m(x,y):m\geq0\}=\{(x,y),p(x,y),p^2(x,y),\ldots\},\] where 
$p^m(x,y)$ denotes the $m$-fold composition of $p$ with itself.  We use Edmonds' permutation technique to show that there is a 2-cell embedding of $\pmn$ (and their special cases $\pn$ and $\bpn$) on a surface where certain cycles correspond to regions. So we obtain a lower bound for the number of regions and then utilizing the Euler-Poincar\'{e} formula, we obtain an upper bound for $\gamma(\pmn)$. We analyze the cases $\gamma(\pn)$ and $\gamma(\bpn)$ separately due to the relevance of these two graphs.

%\hl{[CONSISTENT LABEL VS LABELING?]}

\section{Labeling algorithm for \texorpdfstring{$\pn$}{Pn} and \texorpdfstring{$\bpn$}{BPn}}\label{sec:labelingpnbpn}

In this section, we describe algorithms to label the vertices of the pancake graph $\pn$ and the burnt pancake graph $\bpn$. The labeling produced has the property that the vertices of certain cycles alternate labels. Later, we will find appropriate rotation systems so that these cycles that alternate labels become regions of a 2-cell embedding.% To avoid dealing with the trivial cases $\mathbb{P}_1\cong \mathbb{P}_1(1),\mathbb{P}_2\cong \mathbb{P}_1(2)$, and $\mathbb{BP}_1\cong \mathbb{P}_2(1)$ we assume that $n\geq 3$ if $m=1$, and that $n\geq2$ if $m\geq2$.

\subsection{Base cycles}

In the authors' study of the cycle structure of generalized pancake graphs~\cite{BBP19, BB23}, there was the need for cycles, usually referred to as \emph{base cycles} that traverse all the copies of the $n$ copies of $\mathbb{P}_{n-1}$ contained in $\pn$ and all the $2n$ copies of $\mathbb{BP}_{n-1}$ contained in $\bpn$, owing to the recursive structure of both graphs.

\begin{rem}
 As a matter of notation, we will denote the paths in our graphs by the label of the edges traversed. For example, if $v\in S(m,n)$, $vr_1r_3r_2$ refers to the path that starts at $v$, then goes to  vertex $r_1(v)$, then goes to  vertex $r_3(r_1(v))$, and ends at vertex $r_2(r_3(r_1(v)))$. If the vertex $v$ is omitted, then it may be presumed to be the identity. As we mentioned in the introduction, Section~\ref{sec:prelim-cycles}, we denote cycles by the edges traversed in them. 
\end{rem} 

The following proposition describes the base cycles that will be needed in our algorithms. 

\begin{prop}\label{prop:basecycles} The following are base cycles of $\pn$ and $\bpn$. 
\begin{enumerate}[(i)]
%\item (\cite[Lemma 3.2]{BB23}) If $n>2$, $(r_{n-1}\flop_n)^{mn}$ is a cycle of length $2mn$ in $\pmn$.
\item(\cite[Theorem 3.1.4(c)]{BBrelations}) If $n>2$, 
    $(r_{n-1}r_{n})^n$ is a cycle of length $2n$ in $\pn$.
\item (\cite[Theorem 4.1.5]{BBrelations}) If $n>1$, $(r_{n-1}r_{n})^{2n}$ is a cycle of length $4n$ in $\bpn$.
 %That is, if we traverse the cycle $(r_{n-1}r_{n})^n$ starting from $e_n$, no vertex will be repeated except for the last vertex.
\end{enumerate} In particular, notice that $(r_{n-1}r_{n})^n$ traverses exactly one edge in each of the $n$ copies of $\mathbb{P}_{n-1}$ contained in $\pn$. Similarly, $(r_{n-1}r_{n})^{2n}$ traverses exactly one edge in each of the $2n$ copies of $\mathbb{BP}_{n-1}$ contained in $\mathbb{BP}_n$.
\end{prop}

\subsection{Base cases of our algorithm}

We note that the authors in~\cite{BB22} describe certain cycles in $\pn$ and $\bpn$, in the context of algebraic relations on the set of prefix reversals. In particular, in $\pn$, each path traversed by $(r_2r_3)^3$ is a 6-cycle ~\cite[Theorem 3.1(R2)]{BB22} and that in $\bpn$, each path traversed by $(r_1r_2)^4$ is an 8-cycle~\cite[Theorem 4.1(Rb3) for $k=2$]{BB22}.

It is useful in the discussion of the labeling algorithm to have the following terminology.
\begin{defn} Let $G$ be a graph and let $C$ be a cycle in $G$. We say that a labeling $L$ of the vertices in $G$ is \emph{$C$-alternating} if for every instance of $C$ in $G$, two adjacent vertices of $C$ have different labels in $L$. Furthermore, if $C$ is of the form $(e_1e_2)^k$, for some $k$, we say that $L$ is $\{e_1,e_2\}$-alternating.  
\end{defn}

%We are going to define labelings of the vertices of our graphs that satisfy the following properties, which we name \textit{Property A}.

%\begin{defn}
%We say that a labeling $L$ of $V(G)$ satisfies \emph{Property A} if one of the following is true.
%\begin{itemize}
%    \item Let $G=\pn$ and $C$ be any cycle $\{r_2,r_3\}$-cycle in $\pn$. Then $L$ is $C$-alternating. Due to relation (R3), all such cycles $C$ are 6-cycles. 
%    \item Let $G=\bpn$ and $C$ be any cycle traversing edges that alternate between $r_1$ and $r_2$ in $\bpn$. Then $L$ is $C$-alternating. Due to relation (Rb3), all such cycles $C$ are 8-cycles. 
%    \item Let $G=\pmn$ and $C$ be any cycle traversing edges that alternate between $r_1$ and $r_2$ in $\pmn$. Then $L$ is $C$-alternating. The length of these cycles $C$ is given in Proposition~\ref{prop:r1r2cyclespmn}.
%\end{itemize}    
%\end{defn}

Our goal is to obtain an $\{r_2,r_3\}$-alternating labeling for $V(\pn)$ with $n\geq3$. Additionally, we want to obtain an $\{r_1,r_2\}$-alternating labeling for $V(\bpn)$, with $n\geq2$. To find such labelings, we define the following sets. 
\begin{align*}
    L^1_{even,n}&=\{(r_2r_3)^i,(r_3r_2)^i:0\leq i\leq 1\},\\
    L^1_{odd,n}&=\{(r_2r_3)^ir_2,(r_3r_2)^ir_3:0\leq i\leq 1\},\\
    L^2_{even,n}&=\{(r_1r_2)^i,(r_2r_3)^i:0\leq i\leq2\},\text{ and}\\
    L^2_{odd,n}&=\{(r_2r_1)^ir_2,(r_1r_2)^ir_1:0\leq i\leq 1\}  
    % L^{1,m}_{even,n}&=\{(r_1r_2)^i,(\flop_2\flop_1)^i:0\leq i\leq j\},&\text{ if }m=3j, j\in\mathbb{Z}^+,\\
    % L^{1,m}_{odd,n}&=\{(r_1r_2)^ir_1,(\flop_2\flop_1)^i\flop_2:0\leq i\leq j-1\},&\text{ if }m=3j, j\in\mathbb{Z}^+,\\
    % L^{2,m}_{even,n}&=\{(r_1r_2)^i,(\flop_2\flop_1)^i:0\leq i\leq m\},&\text{ if }m>3,m\neq 3j, j\in\mathbb{Z}^+,\text{ and}\\
    % L^{2,m}_{odd,n}&=\{(r_1r_2)^ir_1,(\flop_2\flop_1)^i\flop_2:0\leq i\leq m-1\},&\text{ if }m>3, m\neq 3j, j\in\mathbb{Z}^+\\
\end{align*}

The intuition behind the definition of these sets is to traverse cycles and have those that are at an even distance away from a base vertex $v$ to be given one label (we use $V_1$) and those vertices in the cycles at an odd distance away from $v$ to be given a different label (we use $V_2$). 

%We define our labeling of $\pn$ recursively. Intuitively, all the $(r_2r_3)^3$ cycles become regions in an embedding. So we want all edges in cycles that have this in $\pn$ to alternate labels. The base of induction is the cycle 
%
%\[\{123\cdot b_n, 213\cdot b_n, 321\cdot b_n,312\cdot b_n,231\cdot b_n,321\cdot b_n\},\]
%
%where $b_n$ is the string $456\cdots n$. Our labeling scheme will maintain two sets $V_1,V_2$, one corresponding to each label, that will partition $\sn$, and in the base of induction,
%\begin{align*}
%    V_1&=\{123\cdot b_n,312\cdot b_n,231\cdot b_n\}\\
%    V_2&=\{213\cdot b_n, 321\cdot b_n,321\cdot b_n\} 
%\end{align*}
%
%In general, we will label any 6-cycle of the form $(r_2r_3)^3$ utilizing a ``base vertex" $v$ by simply alternating the labels starting from $v$. We formalize the labeling procedure procedure in Algorithm~\ref{alg:6l}. 

\begin{algorithm}[H]
\caption{\textsc{AlCyc}$(G,C,v)$, standing for`` \textbf{al}ternating \textbf{cyc}lic" labeling, produces a labeling of the unique $\{r_2,r_3\}$-cycle containing $v$ if $G=\pn$ with $n\geq3$ and the unique $\{r_1,r_2\}$-cycle containing $v$ if $G=\bpn$ with $n\geq 2$.}
\label{alg:base-cases}
\begin{algorithmic}
\Require $G\in\{\pn,\bpn\}$ 
\Require $C$ has the form $(r_2r_3)^3$ if $G=\pn$ , or the form $(r_1r_2)^4$ if $G=\bpn$.%, or given by $(r_1r_2)^{\ell}$ if $G=\pmn$ with $m>2$ and $\ell$ is as in Proposition~\ref{prop:r1r2cyclespmn}
\Require $v\in C$
\State $V_1\gets \emptyset$
\State $V_2\gets \emptyset$
\If {$G=\pn$}
    \State $V_1\gets\{v\cdot w:w\in L^1_{even, n}\}$
    \State $V_2\gets\{v\cdot w:w\in L^1_{odd, n}\}$
\ElsIf{$G=\bpn$}
    \State $V_1\gets\{v\cdot w:w\in L^2_{even, n}\}$
    \State $V_2\gets\{v\cdot w:w\in L^2_{odd, n}\}$
% \ElsIf{$G=\pmn, m>2$, $m=3j$ for some $j\in\mathbb{Z}^+$}
%     \State $V_1\gets\{v\cdot w:w\in L^{1,m}_{even, n}\}$
%     \State $V_2\gets\{v\cdot w:w\in L^{1,m}_{odd, n}\}$
% \ElsIf{$G=\pmn, m>2$, $m\neq 3j$ for some $j\in\mathbb{Z}^+$}
%     \State $V_1\gets\{v\cdot w:w\in L^{2,m}_{even, n}\}$
%     \State $V_2\gets\{v\cdot w:w\in L^{2,m}_{odd, n}\}$
\EndIf
\State\Return $(V_1,V_2)$
\end{algorithmic}
\end{algorithm}

\begin{ex}
We provide some examples of the output of $\AlCyc$.
\begin{description}
    \item[$\AlCyc(\mathbb{P}_3,(r_2r_3)^3,e)$] yields 
    \begin{align*}
        V_1&=\{e,r_2r_3,r_3r_2\},\text{ and}\\
        V_2&=\{r_2,r_3, r_2r_3r_2=r_3r_2r_3\}.
    \end{align*}
    \item[$\AlCyc(\mathbb{BP}_2,(r_1r_2)^4,e)$] yields
        \begin{align*}
        V_1&=\{e,r_1r_2,r_2r_1,(r_1r_2)^2=(r_2r_1)^2\},\text{ and}\\
        V_2&=\{r_1,r_2,r_1r_2r_1,r_2r_1r_2\}.
    \end{align*}
    % \item[$\AlCyc(\mathbb{P}_3(2),(r_1r_2)^2,e)$] yields 
    %     \begin{align*}
    %     V_1&=\{e,r_1r_2=\flop_2\flop_1\},\text{ and}\\
    %     V_2&=\{r_1,\flop_2\}.
    % \end{align*}
    % \item[$\AlCyc(\mathbb{P}_4(2),(r_2r_3)^8,e)$] yields 
    %     \begin{align*}
    %     V_1&=\{e,r_1r_2,\flop_2\flop_1,(r_1r_2)^2,(\flop_2\flop_1)^2,(r_1r_2)^3,(\flop_2\flop_1)^3,(r_1r_2)^4=(\flop_2\flop_1)^4\},\text{ and}\\
    %     V_2&=\{r_1,\flop_2,r_1r_2r_1,\flop_2\flop_1\flop_2,(r_1r_2)^2r_1,(\flop_2\flop_1)^2\flop_2,(r_1r_2)^3r_1,(\flop_2\flop_1)^3\flop_2\}.
    % \end{align*}
\end{description}
\end{ex}

\subsection{Constructor algorithm}

We are at last ready to formalize the labeling procedure to generate the labels that will be used in our rotation systems. In the description of the algorithm that follows, we use $\cup'$ to denote the following tuple operator: If $(V_1,V_2)$ and $(V_3,V_4)$ are two tuples, then \[(V_1,V_2)\cup'(V_3,V_4)=(V_1\cup V_3, V_2\cup V_4).\]

We use $\AlCyc$ to label all the instances of the $\{r_2,r_3\}$-alternating cycles in $\pn$, with $n\geq3$, and all the $\{r_1,r_2\}$-alternating cycles in $\bpn$, with $n\geq2$. We formalize the procedure in the following algorithm. 

\begin{algorithm}[H]
\caption{$\AlGra(G,v,n)$, standing for ``\textbf{al}ternating \textbf{gra}ph" labeling, produces a labeling of $V(G)$ starting with a base vertex $v$}
\label{alg:labels}
\begin{algorithmic}
\Require $G=\pn$ with $n\geq3$ or $G=\bpn$ with $n\geq2$
\Require $v\in G$
\State $V_1\gets\emptyset$
\State $V_2\gets\emptyset$
\If{$G=\pn$}
    \If{$n=3,v=v_1\cdot s_{n-3}, G=\mathbb{P}_3\cdot s_{n-3}$, where $s_{n-3}$ is a suffix of length $n-3$}
        \State $(V_1,V_2)\gets\AlCyc(G,G,v)$
    \ElsIf{$n>3$}
        \For{$i$ from 0 to $n-1$}
            \State $v\gets (r_{n-1}r_n)^i$
            \State $(V_1,V_2)\gets (V_1,V_2) \cup' \AlGra(\mathbb{P}_{n-1}\cdot a,v'\cdot a,n-1$), where $v=v'\cdot a$ 
        \EndFor
    \EndIf
\EndIf
\If{$G=\bpn$}
    \If{$n=2,v=v_1\cdot s_{n-2}, G=\mathbb{BP}_2\cdot s_{n-2}$, where $s_{n-2}$ is a suffix of length $n-2$}
        \State $(V_1,V_2)\gets\AlCyc(G,G,v)$
    \ElsIf{$n>2$}
        \For{$i$ from 0 to $2n-1$}
            \State $v\gets (r_{n-1}r_n)^i$
            \State $(V_1,V_2)\gets (V_1,V_2) \cup' \AlGra(\mathbb{BP}_{n-1}\cdot a,v'\cdot a,n-1$), where $v=v'\cdot a$ 
        \EndFor
    \EndIf
\EndIf
% \If{$G=\pmn, m\geq2$ (this includes $\bpn$)}
%     \If{$n=2,v=v_1\cdot s_{n-2}, G=\mathbb{P}_m(2)\cdot s_{n-2}$, where $s_{n-2}$ is a suffix of length $n-2$}
%         \State $V_1,V_2\gets\AlCyc(G,G,v)$
%     \ElsIf{$n>3$}
%         \For{$i$ from 0 to $mn-1$}
%             \State $v\gets (r_{n-1}r^{-1}_n)^i$
%             \State $V_1,V_2\gets V_1,V_2 \cup' \AlGra(\mathbb{P}_m(n-1)\cdot a,v'\cdot a,n-1$), where $v=v'\cdot a$ 
%         \EndFor
%     \EndIf
%\EndIf
\State\Return $(V_1,V_2)$
\end{algorithmic}
\end{algorithm}

To label the vertices of the graphs we are interested in, we run $\AlGra(\pn, e,n)$, $\AlGra(\bpn, e,n)$ where, by abuse of notation, we denote the identity in each of the groups $S_n$ and $B_n$ by $e$.

We illustrate our procedure by running $\AlGra(\mathbb{P}_4, 1234,4)$ in the following example.

\begin{ex}
We run $\AlGra(\mathbb{P}_4,1234,4)$. The algorithm's first call will initialize $V_1,V_2=\emptyset$ and in the recursive step, it will call $\AlGra(\mathbb{P}_3\cdot 4,123\cdot 4,3)$, which is a 6-cycle, so it will call $\AlCyc(\mathbb{P}_3\cdot 4,123\cdot 4)$ to label $\mathbb{P}_3\cdot 4$. Then, at the second call of the \textbf{for} loop, $v=r_3r_4=4123$, so the algorithm will call $\AlGra(\mathbb{P}_3\cdot 3,412\cdot 3,3)$ which will call $\AlCyc(\mathbb{P}_3\cdot 3,412\cdot 3)$ and label the corresponding 6-cycle $\mathbb{P}_3\cdot 3$. Next, $v=(r_3r_4)^2=3412$ and $\AlGra(\mathbb{P}_3\cdot 2,341\cdot 2,3)$ is called, thus calling $\AlCyc(\mathbb{P}_3\cdot 2,341\cdot 2)$ and labeling the corresponding 6-cycle $\mathbb{P}_3\cdot 2$. In the last run of the \textbf{for} loop, $v=(r_3r_4)^3=2341$ and in the last recursive call, $\AlCyc(\mathbb{P}_3\cdot1,234\cdot 1)$ will be executed, thus labeling the 6-cycle $\mathbb{P}_3\cdot1$. Now the entirety of the vertices of $\mathbb{P}_4$ are labeled. The final labeling is shown in Figure~\ref{f:P4}. Notice that each edge in a cycle of the form $(r_2r_3)^3$ has its endpoints with different labels. 
\end{ex}

\begin{figure}[ht!]
    \begin{center}
    \begin{tikzpicture}[scale=1.5,every node/.style={scale=1}]
        \node (e) at (-1.5,3.5) {$\boxed{1234}$};
            
        \node (2) at (-2.5,2.5) {$\circled{3214}$};
        \node (1) at (-0.5,2.5) {$\circled{2134}$};
        \node (3) at (1.5,3.5) {$\circled{4321}$};
            
        \node (32) at (2.5,-2.5) {$\boxed{4123}$};
        \node (12) at (-2.5,1.5) {$\boxed{2314}$};
        \node (21) at (-0.5,1.5) {$\boxed{3124}$};
        \node (31) at (-0.5,-2.5) {$\circled{4312}$};
        \node (13) at (0.5,2.5) {$\boxed{3421}$};
        \node (23) at (2.5,2.5) {$\boxed{2341}$};

        \node (232) at (1.5,-3.5) {$\circled{2143}$};						
        \node (132) at (2.5,-1.5) {$\circled{1423}$};
        \node (312) at (-2.5,-1.5) {$\boxed{4132}$};
        \node (212) at (-1.5,0.5) {$\circled{1324}$};
        \node (321) at (0.5,-1.5) {$\circled{4213}$};
        \node (231) at (-0.5,-1.5) {$\boxed{1342}$};
        \node (131) at (-1.5,-3.5) {$\boxed{3412}$};
        \node (213) at (0.5,1.5) {$\circled{2431}$};
        \node (313) at (0.5,-2.5) {$\boxed{1243}$};
        \node (123) at (2.5,1.5) {$\circled{3241}$};
        \node (323) at (-2.5,-2.5) {$\circled{1432}$};
            
        \node (1321) at (1.5,-0.5) {$\boxed{2413}$};
        \node (1231) at (-1.5,-0.5) {$\circled{3142}$};
        \node (1213) at (1.5,0.5) {$\boxed{4231}$};
            
        \draw[line width=2.pt,red] (e)--(1) (2)--(12) (3)--(13) (31)--(131) (21)--(212) (32)--(132) (23)--(123) (231)--(1231) (321)--(1321) (213)--(1213);
        \draw[line width=2.pt,blue]  (1)--(21)  (31)--(231) (12)--(212)  (13)--(213) (312)--(1231) (132)--(1321) (123)--(1213);
        \draw[line width=2.pt,purple!50!black]   (21)--(321) (12)--(312)  (212)--(1213);% (1) arc (100:260:2 and 2.5) (313) arc (-80:80:2 and 2.5); %(1)--(31) (13)--(313);
        \draw[line width=2.pt,purple!50!black,domain=-72:72] plot ({0.25+2*cos(\x)},{2.5*sin(\x)}) plot ({-0.25+2*cos(\x+180)},{2.5*sin(\x+180)});

        \draw[line width=2.pt] (232) edge[red] (313) 
                  (232) edge[purple!50!black,dashed] (131) 
                  (132) edge[purple!50!black] (123) 
                  (321) edge[blue] (313) 
                  (231) edge[purple!50!black] (213)
                  (312) edge[red] (323)
                  (131) edge[blue,dashed] (323)
                  (1321) edge[purple!50!black] (1231);

        \draw[line width=2.pt,blue,dashed] (e)--(2) (32)--(232) (3)--(23);

        \draw[line width=2.pt,purple!50!black,dashed] (e)--(3) (2)--(32) (23)--(323);

        %\draw[line width=2.pt,purple!50!black,domain=-72:72,dashed]
    \end{tikzpicture}	
    \caption{Labeling of $V(\mathbb{P}_4)$ after running $\AlGra(\mathbb{P}_4,1234, 4)$. The vertices in $V_1$ are indicated with a rectangle and the vertices in $V_2$ are indicated with a rounded rectangle. The base cycle $(r_3r_4)^4$ starting at $1234$ is indicated with dashed edges.}
    \label{f:P4}
    \end{center}
\end{figure}

\begin{figure}[ht!]
    \centering
    \begin{tikzpicture}[scale=0.19,line cap=round,line join=round,>=triangle 45,x=1.0cm,y=1.0cm]
    	\begin{scriptsize}
    	% ending in 3
    	\node (e) at (-13,-33) {$\boxed{123}$};
    	\node (0) at (13,-33) {$\circled{\underline{1}23}$};
    	\node (10) at (33,-14) {$\boxed{\underline{2}13}$};
    	\node (010) at (33,14) {$\circled{213}$};
    	\node (1010) at (13,33) {$\boxed{\underline{1}\,\underline{2}3}$};
    	\node (101) at (-13,33) {$\circled{1\underline{2}3}$};
    	\node (01) at (-33,14) {$\boxed{2\underline{1}3}$};
    	\node (1) at  (-33,-14) {$\circled{\underline{2}\,\underline{1}3}$};
    
    	%ending in -3
    	\node (202) at  (-4,-9) {$\boxed{12\underline{3}}$};
    	\node (0202) at  (4,-9) {$\circled{\underline{1}2\underline{3}}$};
    	\node (10202) at  (9,-3.5) {$\boxed{\underline{2}1\underline{3}}$};
    	\node (010202) at  (9,4) {$\circled{21\underline{3}}$};
    	\node (1010202) at  (4,9) {$\boxed{\underline{1}\;\underline{2}\,\underline{3}}$};
    	\node (101202) at  (-4,9) {$\circled{1\underline{2}\,\underline{3}}$};
    	\node (01202) at  (-9,4) {$\boxed{2\underline{1}\;\underline{3}}$};
    	\node (1202) at  (-9,-3.5) {$\circled{\underline{2}\,\underline{1}\,\underline{3}}$};
    
    	%ending in 2
    	\node (21) at  (-12.5,-20) {$\boxed{\underline{3}12}$};
    	\node (121) at  (-6,-15) {$\circled{\underline{1}32}$};
    	\node (0121) at  (6,-15) {$\boxed{132}$};
    	\node (10121) at  (12.5,-20) {$\circled{\underline{3}\,\underline{1}\,2}$};
    	\node (010121) at  (11,-25) {$\boxed{3\underline{1}2}$};
    	\node (01021) at  (7,-29.5) {$\circled{1\underline{3}2}$};
    	\node (1021) at  (-7,-29.5) {$\boxed{\underline{1}\,\underline{3}2}$};
    	\node (021) at  (-11,-25) {$\circled{312}$};
    	
    	%ending in -2
    	\node (212) at (6,14.75) {$\circled{1\underline{3}\,\underline{2}}$};
    	\node (0212) at (-6,14.75) {$\boxed{\underline{1}\,\underline{3}\,\underline{2}}$};
    	\node (10212) at (-12.5,20) {$\circled{31\underline{2}}$};
    	\node (010212) at (-12,25) {$\boxed{\underline{3}1\underline{2}}$};
    	\node (1010212) at (-7,29.5) {$\circled{\underline{1}3\underline{2}}$};
    	\node (101212) at (7,29.5) {$\boxed{13\underline{2}}$};
    	\node (01212) at (12,25) {$\circled{\underline{3}\,\underline{1}\,\underline{2}}$};
    	\node (1212) at (12.5,20) {$\boxed{3\underline{1}\,\underline{2}}$};
    	
    	%ending in 1
    	\node (20) at  (29.5,-8) {$\circled{\underline{3}\,\underline{2}1}$};
    	\node (020) at  (29.5,7) {$\boxed{3\underline{2}1}$};
    	\node (1020) at  (26,11.5) {$\circled{2\underline{3}1}$};%
    	\node (01020) at  (18.5,11.5) {$\boxed{\underline{2}\,\underline{3}1}$};
    	\node (101020) at  (17,5.5) {$\circled{321}$};
    	\node (10120) at  (16.5,-7) {$\boxed{\underline{3}\,21}$};	
    	\node (0120) at  (18.5,-11.5) {$\circled{\underline{2}\,31}$};
    	\node (120) at  (24.5,-11.5) {$\boxed{231}$};	
    
    	%ending in -1
    	\node (2) at  (-29.5,-8) {$\circled{\underline{3}\,\underline{2}\,\underline{1}}$};
    	\node (02) at  (-29.5,7) {$\boxed{3\underline{2}\,\underline{1}}$};
    	\node (102) at  (-26,11.5) {$\circled{2\underline{3}\,\underline{1}}$};%
    	\node (0102) at  (-18.5,11.5) {$\boxed{\underline{2}\,\underline{3}\,\underline{1}}$};
    	\node (10102) at  (-17,5.5) {$\circled{32\underline{1}}$};
    	\node (1012) at  (-16.5,-7) {$\boxed{\underline{3}2\underline{1}}$}; 
    	\node (012) at  (-18.5,-11.5) {$\circled{\underline{2}3\underline{1}}$};
    	\node (12) at  (-24.5,-11.5) {$\boxed{23\underline{1}}$};	
    	\end{scriptsize}
    	
    	% f_0
    	\draw[line width=2,purple!50!black] (e)--(0) (10)--(010) (101)--(1010) (1)--(01) (202)--(0202) (10202)--(010202) 
    	(101202)--(1010202) (1202)--(01202) (21)--(021) (121)--(0121) (10121)--(010121) (1021)--(01021) 
    	(212)--(0212) (10212)--(010212) (101212)--(1010212) (1212)--(01212) (20)--(020) (1020)--(01020) 
    	(10120)--(101020) (120)--(0120) (2)--(02) (102)--(0102) (1012)--(10102) (12)--(012);
    	% f_1
    	\draw[line width=2,red] (0)--(10) (01)--(101) (02)--(102) (20)--(120) 
    	(010)--(1010) (012)--(1012) (020)--(1020) (021)--(1021) (010212)--(1010212) 
    	(0121)--(10121) (202)--(1202) (0202)--(10202) (01202)--(101202) 
    	(0102)--(10102) (01021)--(010121) (0120)--(10120) (01212)--(101212) (0212)--(10212);
    	% f_2
    	\draw[line width=2,blue] (0)--(20) (10)--(10121) (02)--(202) (010)--(01212) 
    	(1010)--(10120) (101)--(1012) (01)--(010212) (020)--(0202)  (010121)--(10202) (10102)--(101202)
    	(120)--(0212) (0121)--(0102)
    	(01021)--(012) (1021)--(0120) (021)--(1202) (101212)--(102) (1010212)--(1020) (01202)--(10212);

            % base cycle
            \draw[line width=2,red,dashed] (e)--(1) (21)--(121) (01020)--(101020) (010202)--(1010202) (212)--(1212) (2)--(12);
            
            \draw[line width=2,blue,dashed] (1)--(21) (121)--(01020) (101020)--(1010202) (010202)--(1212) (12)--(212) (e)--(2);

    \end{tikzpicture}
    \caption{Labeling of $V(\mathbb{BP}_4)$ after running $\AlGra(\mathbb{BP}_3,123,3)$. The vertices in $V_1$ are indicated with a rectangle and the vertices in $V_2$ are indicated with a rounded rectangle. The base cycle $(r_2r_3)^6$ starting at $123$ is indicated with dashed edges.}
    \label{f:BP3}
\end{figure}
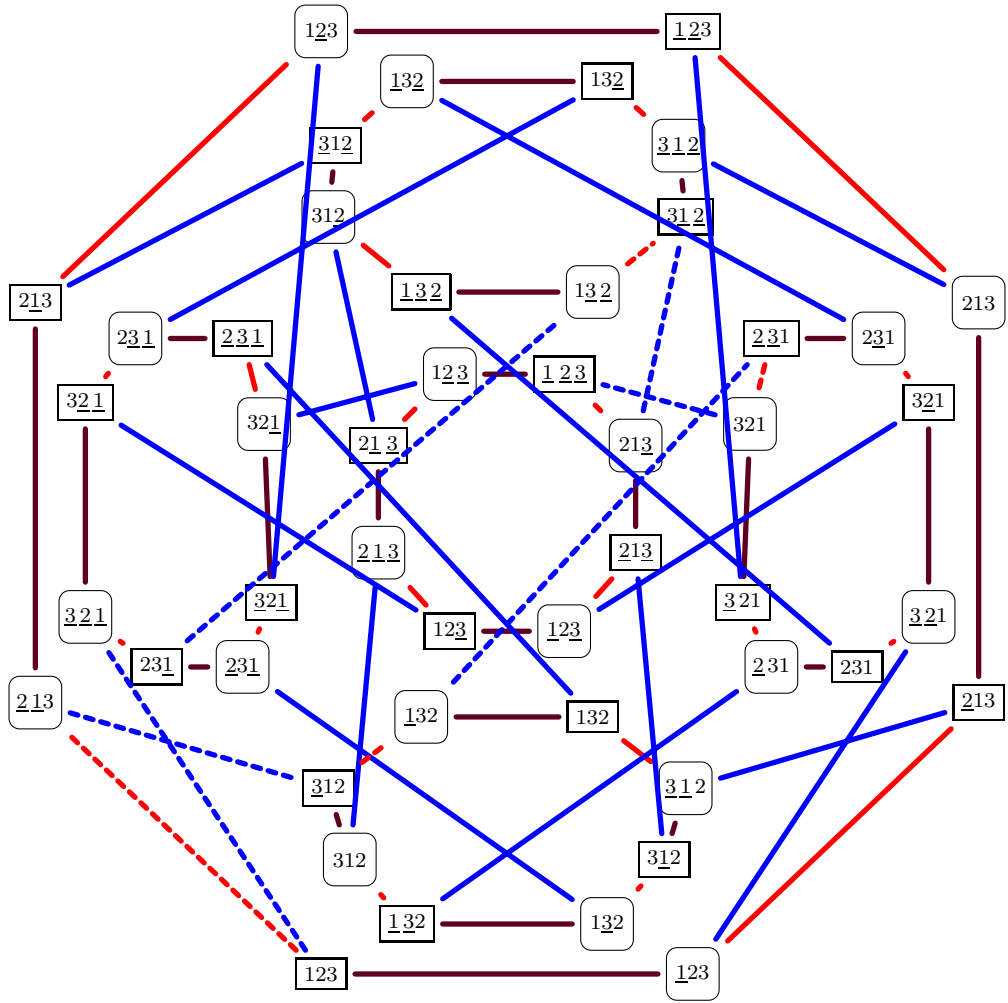

We also present the end results of our labeling procedure $\AlGra$ on $\mathbb{BP}_3$ in Figure~\ref{f:BP3}. 

We now prove that $\AlGra$ produces a $C$-alternating labeling, where $C$ is any $\{r_2,r_3\}$-cycle in the case of $\pn$ or a $\{r_1,r_2\}$-cycle in the case of $\bpn$.

\begin{prop} $\AlGra(G,e,n)$, with $G\in\{\pn,\bpn\}$, labels all of the vertices of $G$. Moreover, 
\begin{enumerate}
    \item[(i)] If $G=\pn$, with $n\geq3$, then $\AlGra$ produces an $\{r_2,r_3\}$-alternating labeling.
    \item[(ii)] If $G=\bpn$ with $n\geq2$, then $\AlGra$ produces an $\{r_1,r_2\}$-alternating labeling.
\end{enumerate}
\end{prop}
\begin{proof} We proceed by induction to show that $\AlGra$ labels all the vertices of $G$. In the base case, $G$ is isomorphic to either $\mathbb{P}_3$ or $\mathbb{BP}_2$. These two graphs are a 6-cycle or an 8-cycle, respectively, and are labeled by $\AlCyc$ in such a way that the labels of the cycles alternate between $V_1$ and $V_2$. As inductive hypothesis, we assume that $\AlGra$ properly labels all copies of $\mathbb{P}_{n-1}$ and $\mathbb{BP}_{n-1}$ that are contained in $\mathbb{P}_n$ or $\mathbb{BP}_n$, respectively. Notice that the \textbf{for} loops of the algorithm will change the base vertex $v$ to take the form $(r_{n-1}r_n)^i$. As $i$ varies, $v$ will become a vertex in each of the different copies of $\mathbb{P}_{n-1}$ and $\mathbb{BP}_{n-1}$ contained in $\mathbb{P}_n$ and $\mathbb{BP}_n$, respectively. Therefore, via the induction hypothesis, all of the vertices in $\mathbb{P}_n$ and $\mathbb{BP}_n$ are labeled since $v$ will become the base vertex of each of the copies of $\mathbb{P}_{n-1}$ and $\mathbb{BP}_{n-1}$ and $\AlGra$ labels all of the vertices of each of those copies contained in $\pn$ or $\bpn$, respectively. 

    For (i), we also proceed by induction. The base case of $\AlGra$ will call upon $\AlCyc$ and thus one obtains an $\{r_2,r_3\}$-alternating labeling of a 6-cycle isomorphic to $\mathbb{P}_3$. We suppose as inductive hypothesis that $\AlGra$ produces an $\{r_2,r_3\}$-alternating labeling of each of the $\mathbb{P}_{n-1}$ copies contained in $\pn$. The inductive step is established by noticing that the \textbf{for} loop will change the base vertex $v$ to be a vertex in each of the copies of $\mathbb{P}_{n-1}$ inside of $\pn$. Since each of $\mathbb{P}_{n-1}$ copies has been given an $\{r_2,r_3\}$-alternating labeling, the entire graph $\pn$ has an $\{r_2,r_3\}$-alternating labeling. 
    
The argument for (ii) is the same, \emph{mutatis mutandis}, for $\{r_1,r_2\}$-cycles in $\bpn$.
\end{proof}

% \begin{prop}[Theorem 4.1.5 \cite{BBrelations}]
%     $(r_{n-1}r_{n})^{2n}$ is a cycle of length $4n$ in $\bpn$. That is, if we traverse the cycle $(r_{n-1}r_{n})^{2n}$ starting from $e_n$, no vertex will be repeated except for the last vertex. 
% \end{prop}

% \subsection{Labeling $\bpn$}

% \begin{algorithm}[H]
% \caption{$8L(C,v)$ to label a eight-cycle $C$ of the form $(r_1r_2)^4$ in $\bpn$ given, a ``base vertex" $v$ in $C$.}
% \label{alg:8l}
% \begin{algorithmic}
% %\State $s_n\gets $ largest suffix shared by the endpoints of $b$.
% \State $V_1\gets\{v,vr_2r_1,vr_1r_2,v(r_1r_2)^2=v(r_2r_1)^2\}$
% \State $V_2\gets\{vr_1,vr_2,vr_1r_2r_1,vr_2r_1r_2\}$
% \State \Return $V_1,V_2$
% \end{algorithmic}
% \end{algorithm}

% After the base of induction, the algorithm will call upon itself to label each of the copies of $\mathbb{BP}_{n-1}$ in $\bpn$. We formalize the procedure in Algorithm~\ref{alg:LGB}.

% \begin{algorithm}[H]
% \caption{$LGB(v,G,n)$ to label $\bpn$, $n\geq2$, and a ``base vertex" $v\in V(G)$.}\label{alg:LGB}
% \begin{algorithmic}
% \State $V_1\gets\emptyset$
% \State $V_2\gets\emptyset$
% \If{$n=2, v=v_1\cdot s_{n-2},G=\mathbb{BP}_2\cdot s_{n-2}$ (where $s_{n-2}$ is a suffix of length $n-3$)}
% \State $V_1,V_2\gets 8L(G,v)$
% \EndIf
% \For{$i$ from 0 to $n-1$}
% \State $v\gets (r_{n-1}r_n)^i$
% \State{/*We use $\cup'$ to denote the following tuple operator: $V_1,V_2\cup'V_3,V_4=V_1\cup V_3, V_2\cup V_4$/*}
% \State $V_1,V_2\gets V_1,V_2 \cup' LG(v'\cdot a,\mathbb{BP}_{n-1}\cdot a,n-1$), where $v=v'\cdot a,a\in[n]$  
% \EndFor
% \State\Return $V_1,V_2$
% \end{algorithmic}
% \end{algorithm}

\section{Genus bounds}\label{sec:bounds}

In this section present our bounds for the genus of $\pn,\bpn$, and $\pmn$. We first discuss the genus of $\pn$, then the genus of $\bpn$, and finally the genus of $\pmn$ for $m>2$. The techniques used to establish the bounds for $\gamma(\pn)$ and $\gamma(\bpn)$ are similar in spirit, whereas the proof of the bound for $\gamma(\pmn)$ makes use of certain cycles that are not present if $m\leq2$. We also note that the graphs $\mathbb{P}_1$, $\mathbb{P}_2$, $\mathbb{P}_3$, $\mathbb{BP}_1$, $\mathbb{BP}_2$, and $\mathbb{P}_m(1)$ are all planar and thus have genus 0.

\subsection{Upper bound for \texorpdfstring{$\gamma(\pn$)}{gamma(Pn)}}

The following theorem is our first main result, a tighter upper bound for $\gamma(\pn)$ than what is currently known in the literature, see Nguyen and Bettayeb~\cite{NB09}.

\begin{thm}\label{t:pnupperbound}
    If $n>3$, then
    \[
    \gamma(\pn) \leq n!\left(\frac{3n-10}{12}\right) + 1.
    \]
\end{thm}

\begin{proof}

Our approach is to show that cycles of the form $(r_2r_3)^3$ are orbits under a particular rotation system. Indeed, we consider the following rotation system
\begin{equation}
p_v=
\begin{cases}
    (vr_3, vr_2, vr_4, \ldots)&\text{ if  $v\in V_1$, and}\\
    (vr_2, vr_3, vr_4, \ldots)&\text{ if  $v\in V_2$,}
\end{cases}
\end{equation}
where the sets $V_1$ and $V_2$ are assigned according to Algorithm \ref{alg:labels}.

Now, we trace the orbit of edges of the form $(v,vr_2)$ if $v\in V_1$, which are of the form
\begin{align*}
(v,vr_2) &  \\
p(v,vr_2) &= (vr_2,p_{vr_2}(v))=(vr_2,vr_2r_3)&\text{ since $vr_2\in V_2$}\\
p(vr_2,vr_2r_3) &= (vr_2r_3, p_{vr_2r_3}(vr_2))=(vr_2r_3, vr_2r_3r_2)&\text{ since $vr_2r_3\in V_1$}\\
p(vr_2r_3, vr_2r_3r_2) &= (vr_2r_3r_2,v(r_2r_3)^2)\\
p(vr_2r_3r_2, v(r_2r_3)^2) &= (v(r_2r_3)^2,v(r_2r_3)^2r_2)\\
p(v(r_2r_3)^2,v(r_2r_3)^2r_2) &= (v(r_2r_3)^2r_2,v(r_2r_3)^3)=(vr_3,v)
\end{align*}

As we vary the permutation $v\in V_1$, every vertex in $\pn$ will be contained in an orbit of the above form. Thus, there are $n!/6$ of these orbits, with $n>3$. 

There are certainly more orbits that are not of the above form. So we conclude that the number of regions in the 2-cell embedding given by this rotation system is at least $n!/6$.

The number of vertices of $\pn$ is $n!$ and the number of edges in $\pn$ is $(n-1)n!/2$. Thus, combining the Euler-Poincar\'{e} formula with this lower bound for $r$ gives us
\begin{align*}
    2-2\gamma(\pn) &\geq n! - (n-1)n!/2 + n!/6\\
    &= \frac{1}{6}(10 - 3 n) n!.
\end{align*}
Therefore, solving for $\gamma(\pn)$ we find the upper bound of 
\[
    \gamma(\pn) \leq n!\left(\frac{3n-10}{12}\right) + 1,
\] which is the desired result.
\end{proof}

\begin{rem}\label{rem:pnbounds}
    This is a tighter upper bound than that found in Nguyen and Bettayeb~\cite{NB09}. For reference, these bounds are 
    \[
    n!\left(\frac{n-4}{6}\right) + 1 \leq\gamma(\pn)\leq n!\left(\frac{n-3}{4}\right)-\frac{n}{2}+1.
    \]

    In other words, $\gamma(\pn)$ is $\Theta(nn!)$. Thus, while the bound in Nguyen and Bettayeb~\cite{NB09} is asymptotically tight, our upper bound is $n!/12-n/2$ less than the one above, for $n \geq 5$.
\end{rem}

\begin{rem}
    In the particular case of $\mathbb{P}_4$ and the rotation system described in the proof of Theorem \ref{t:pnupperbound}, i.e., the graph illustrated in Figure \ref{f:P4}, the authors computed the complete set of orbits. With this particular rotation system and labeling there are eight regions. We list the reversals forming the boundary cycles of each region with a seed edge for the orbit.
        \begin{table}[ht!]
            \begin{tabular}{|l|l|}
                \hline
                edge & boundary cycle \\ \hline
                $(1234,2134)$ & $(r_2r_3)^3$\\ \hline
                $(4321,3421)$ & $(r_2r_3)^3$\\ \hline
                $(4312,3412)$ & $(r_2r_3)^3$\\ \hline
                $(4123,1423)$ & $(r_2r_3)^3$\\ \hline
                $(2134,1234)$ & $(r_2r_4r_3r_4r_2r_4)^2$\\ \hline
                $(3124,2134)$ & $(r_3r_4)^4$\\ \hline
                $(1324,3124)$ & $((r_2r_4)^2(r_3r_4)^2)^2$\\ \hline
                $(3214,2314)$ & $(r_2r_4)^4$\\ \hline
            \end{tabular}
        \end{table}
    Although not all possible labelings of the vertices were attempted to find precisely the genus of $\mathbb{P}_4$, there is a strong indication that this is the maximum number of regions yielded from any rotation system. So, we make the following conjecture.
\end{rem}

\begin{conj}
    The graph $\mathbb{P}_4$ is a pretzel graph. That is, $\gamma(\mathbb{P}_4) = 3$.
\end{conj}

In Table~\ref{tab:pn} we summarize some values for the bounds for $\gamma(\pn).$

\begin{table}[ht!]
\begin{tabular}{|c|r|r|r|}
\hline
    & \multicolumn{3}{|c|}{$\gamma(\pn)$} \\ \hline
$n$ & lower bound & new upper bound & $\lfloor$old upper bound$\rfloor$ \\ \hline
4   & 1 & 5 & 5             \\ \hline
5   & 21 & 51 & 58           \\ \hline
6   & 241 & 481 & 538           \\ \hline
7   & 2\;521 & 4\;621 & 5\;037          \\ \hline
8   & 26\;881 & 47\;041 & 50\;397       \\ \hline
9   & 302\;401 & 514\;081 & 544\;317        \\ \hline
10  & 3\;628\;801 & 6\;048\;001 & 6\;350\;396    \\ \hline
11  & 46\;569\;601 & 76\;507\;201 & 79\;833\;595    \\ \hline
12  & 638\;668\;801 & 1\;037\;836\;801 & 1\;077\;753\;595    \\ \hline
\end{tabular}%
\caption{Lower and upper bounds for $\gamma(\pn).$ Notice that our upper bound is tighter than the bound in~\cite{NB09}, for $n\geq5$.}
\label{tab:pn}
\end{table}

\subsection{Lower bound for \texorpdfstring{$\gamma(\bpn)$}{gamma(BPn)}}

A lower bound for the genus of the burnt pancake graphs is easily obtained utilizing known bounds regarding the girth of a graph and the fact that every edge is on the boundary of at least two regions. Recall the following Theorem, see e.g., Chartrand, Lesniak, and Zhang~\cite[Theorem 11.14]{CLZ16}.

\begin{thm}\label{thm:lowerbound}
If $G = (V, E)$ is a connected graph with girth $\alpha$ then
\[\gamma(G)\geq \frac{1}{2}\left(|E|\left(1-\frac{2}{\alpha}\right)-|V|\right)+1
\]
\end{thm}

In the case of $\bpn$, with $n >2$, it is known that $|V(\bpn)|=2^n\cdot n!,|E(\bpn)|=2^{n-1}\cdot n! \cdot n,$ and $\alpha=8$. Thus, utilizing Theorem~\ref{thm:lowerbound} directly, we obtain the following corollary.
\begin{cor}\label{cor:bpnlowerbound}
    $\gamma(\bpn)\geq 2^{n - 4} (3 n - 8) n! + 1$.
\end{cor}

\subsection{Upper bound for \texorpdfstring{$\gamma(\bpn)$}{gamma(BPn)}}

Our second main result is the first known upper bound for $\gamma(\bpn)$, which we give in the following theorem. 

\begin{thm}\label{t:bpnupperbound}
    If $n \geq3$, then
    \[
    \gamma(\bpn) \leq 2^{n-4}(4n-9)n! + 1.
    \]
\end{thm}

\begin{proof}

We proceed by showing that cycles of the form $(r_1r_2)^4$ correspond to orbits under a rotation system. Indeed, we consider the following rotation system
\begin{equation}
p_v=
\begin{cases}
    (vr_1, vr_2, vr_3, \ldots)&\text{ if  $v\in V_1$, and}\\
    (vr_2, vr_1, vr_3, \ldots)&\text{ if  $v\in V_2$},
\end{cases}
\end{equation}
where the sets $V_1$ and $V_2$ are assigned according to Algorithm~\ref{alg:labels}.

We notice that the orbits of edges of the form $(v,vr_2)$ if $v\in V_1$ are of the form
\begin{align*}
(v,vr_2) & \\
p(v,vr_2) &= (vr_2,p_{vr_2}(v))\\
&=(vr_2,vr_2r_1)&\text{ since $vr_2\in V_2$}\\
p(vr_2,vr_2r_1) &= (vr_2r_1, p_{vr_2r_1}(vr_2))\\
&=(vr_2r_1, vr_2r_1r_2)&\text{ since $vr_2r_1\in V_1$}\\
p(vr_2r_1, vr_2r_1r_2) &= (vr_2r_1r_2,v(r_2r_1)^2)\\
p(vr_2r_1r_2, v(r_2r_1)^2) &= (v(r_2r_1)^2,v(r_2r_1)^2r_2)\\
p(v(r_2r_1)^2,v(r_2r_1)^2r_2) &= (v(r_2r_1)^2r_2,v(r_2r_1)^3)\\
p(v(r_2r_1)^2r_2,v(r_2r_1)^3) &= (v(r_2r_1)^3,v(r_2r_1)^3r_2)\\
p(v(r_2r_1)^3,v(r_2r_1)^3r_2) &= (v(r_2r_1)^3r_2,v(r_2r_1)^4)=(vr_1,v)
\end{align*}

As we vary the signed permutation of $v\in V_1$, every vertex in $\bpn$ will be contained in an orbit of the above form. Thus, there are $2^n \cdot n!/8 = 2^{n-3} \cdot n!$ of these orbits. 

Of course, this is not the full extent of orbits in this 2-cell embedding. However, we can bound the number of regions below by the count of these orbits, that is this rotation system yields at least $2^{n-3}\cdot n!$ regions.

The number of vertices of $\bpn$ is $2^{n} \cdot n!$ and the number of edges in $\bpn$ is $2^{n-1}\cdot n! \cdot n$. Thus, combining the Euler-Poincar\'{e} formula with this lower bound for $r$ gives us
\begin{align*}
    2-2\gamma(\bpn) &\geq 2^{n}\cdot n! - 2^{n-1} \cdot n! \cdot n + 2^{n-3}\cdot n!\\
    &= 2^{n-3}(9-4n) n!.
\end{align*}
Therefore, solving for $\gamma(\bpn)$ we find the upper bound of 
\[
    \gamma(\bpn) \leq 2^{n-4}(4n-9)n!+1,
\] as desired.
\end{proof}

\begin{rem}
    In the particular case of $\mathbb{BP}_3$ and the rotation system described in the proof of Theorem \ref{t:bpnupperbound}, i.e., the graph illustrated in Figure \ref{f:BP3}, the authors computed the complete set of orbits. With this particular rotation system and labeling there are ten regions. We list the reversals forming the boundary cycles of each region with a seed edge for the orbit.
        \begin{table}[ht!]
            \begin{tabular}{|l|l|}
                \hline
                edge & boundary cycle \\ \hline
                $(123,\underline{2}\underline{1}3)$ & $(r_1r_2)^4$\\ \hline
                $(\underline{3}12,\underline{1}32)$ & $(r_1r_2)^4$\\ \hline
                $(\underline{2}\underline{3}1,321)$ & $(r_1r_2)^4$\\ \hline
                $(\underline{1}\underline{2}\underline{3},21\underline{3})$ & $(r_1r_2)^4$\\ \hline
                $(3\underline{1}\underline{2},1\underline{3}\underline{2})$ & $(r_1r_2)^4$\\ \hline
                $(23\underline{1},\underline{3}\underline{2}\underline{1})$ & $(r_1r_2)^4$\\ \hline
                $(123,\underline{1}23)$ & $(r_1r_3r_2r_3r_1r_3)^6$\\ \hline
                $(123,\underline{3}\underline{2}\underline{1})$ & $(r_3r_2)^6$\\ \hline
                $(\underline{1}23,\underline{2}13)$ & $((r_2r_3)^2(r_1r_3)^2)^3$\\ \hline
                $(1\underline{2}3, 2\underline{1}3)$ & $(r_2r_3(r_1r_3)^2r_2r_3)^3$\\ \hline
            \end{tabular}
        \end{table}
    Although not all possible labelings of the vertices were attempted to find precisely the genus of $\mathbb{BP}_3$, there is a strong indication that this is the maximum number of regions yielded from any rotation system. So, we make the following conjecture.
\end{rem}

\begin{conj}
    The graph $\mathbb{BP}_3$ is octo-toriodal. That is, $\gamma(\mathbb{BP}_3) = 8$.
\end{conj}

We compute values (see Table~\ref{tab:bpn}) of our lower and upper bounds for $\gamma(\bpn)$ for several values of $n$.

\begin{table}[ht!]
\begin{tabular}{|c|r|r|}
\hline
    & \multicolumn{2}{|c|}{$\gamma(\bpn)$} \\ \hline
$n$ & lower bound & upper bound \\ \hline
3   & 4 & 10            \\ \hline
4   & 97 & 169            \\ \hline
5   & 1\;681 & 2\;641            \\ \hline
6   & 28\;801 & 43\;201           \\ \hline
7   & 524\;161 & 766\;081         \\ \hline
8   & 10\;321\;921 & 14\;837\;761         \\ \hline
9   & 220\;631\;041 & 313\;528\;321       \\ \hline
10  & 5\;109\;350\;401 & 7\;199\;539\;201     \\ \hline
11  & 127\;733\;760\;001 & 178\;827\;264\;001      \\ \hline
12  & 3\;433\;483\;468\;801 & 4\;782\;351\;974\;401     \\ \hline
\end{tabular}%
\caption{Values for lower and upper bounds for $\gamma(\bpn)$ for small values of $n$. }
\label{tab:bpn}
\end{table}

\subsection{Lower bound for \texorpdfstring{$\gamma(\pmn)$}{gamma(Pm(n))},  \texorpdfstring{$m>2$}{m>2}}

In \cite[Theorem 4.2]{BB23}, the authors show the girth of the undirected generalized pancake graph $\pmn$ to be $\min\{6, m\}$. Using this girth we can find a lower bound for the genus of the generalized pancake graphs.

In the case of $\pmn$, with $m\geq 3$ and $n \geq 2$, we have that $|V(\pmn)|=m^n\cdot n!,|E(\pmn)|=m^{n}n! \cdot n,$ and $\alpha=\min\{6,m\}$. Then it follows from Theorem~\ref{thm:lowerbound} that
\begin{cor}\label{cor:pmnlowerbound}
    $\gamma(\pmn)\geq \begin{cases}
        \frac{1}{2}m^{n - 1} ((m-2) n - m) n! + 1, & m \in \{3,4,5\}\\
        \frac{1}{6}m^n (2n-3) n! + 1, & m \geq 6.
    \end{cases}$
\end{cor}

\subsection{Upper bound for \texorpdfstring{$\gamma(\pmn)$}{gamma(Pm(n))}, \texorpdfstring{$m>2$}{m>2}}

To find an upper bound for the generalized pancake graphs, we use the following facts of prefix reversals in $S(m,n)$.

\begin{defn}
    The \emph{order} of a prefix reversal $r_i$ of $S(m,n)$ with $1\leq i\leq n$ is the smallest positive integer $k$ such that $(r^k_i)(e)=e$, where $r^k_i$ denotes the $k$-fold composition of $r_i$ with itself. We use $o(r_i)$ to refer to the order of $r_i$.
\end{defn}

The following is an easily-verified fact about the order of prefix reversals. 

\begin{prop}
    Let $\{r_i\}_{i=1}^n$ be the set of prefix reversals in $S(m,n)$ with $m>2$, $n>1$. Then $o(r_1)=m$ and for $2\leq i\leq n$. Then 
    \[o(r_i)=
    \begin{cases}
        m&\text{ if }m\text{ is even }\\
        2m&\text{ if }m\text{ is odd }
    \end{cases}
    \]
\end{prop}

\begin{rem}
    We use the following notation when traversing the graph $\pmn$.
    \begin{enumerate}[(i)]
        \item If we traverse an edge in $\pmn$ identified with $r_i$ in the opposite direction as it appears in $P(m,n)$, we denote it by $\flop_i$ 
        \item In $\pmn$ if we traverse $k$ edges that have the same label, say $r_i$, then we denote this by $r^k_i$.
    \end{enumerate}
\end{rem}

%in a similar manner as $\pn$ and $\bpn$ above, we will establish the order of cycles using only one of $\{r_1, \flop_1\}$ and one of $\{r_2, \flop_2\}$.

We are now ready to establish an upper bound for $\gamma(\pmn)$.

\begin{thm}\label{t:pmnupperbound}
    If $m \geq 3$ and $n \geq 2$, then
\[\gamma(\pmn)\leq
\begin{cases}
    \frac{1}{2}m^{n-1} (mn-m-n)n!+1,&\text{ if }m\text{ is even}\\
    \frac{1}{2}m^{n-1} (2 m n - 2 m - n-1) n!+1,&\text{ if }m\text{ is odd.}
\end{cases}\] 
\end{thm}
 \begin{proof}
     Consider the following rotation system for all $v\in S(m,n)$.

\begin{equation}\label{eq:rotpmn}
p_v=(v\flop_1,vr_1,v\flop_2,vr_2,\ldots,v\flop_n,vr_n)    
\end{equation}

Direct computation shows that the orbit of $(v,vr_i)$ has the same number of elements as $o(r_i)$. Indeed, the orbit of $(v,vr_i)$ contains
\begin{align*}
    (v,vr_i), &\\
    p(v,vr_i) &= (vr_i,p_{vr_i}(v))=(vr_i,vr^2_i)&\text{since $p_{vr_i}(v)=p_{vr_i}(vr_i\flop_i)$,}\\
    p(vr_i,vr^2_i) &= (vr^2_i,p_{vr^2_i}(vr_i))=(vr^2_i,vr^3_i)&\text{since $p_{vr^2_i}(vr_i)=p_{vr^2_i}(vr^2_i\flop_i)$,}
\end{align*} 
and so on.

Following the computation, the last edge contained in the orbit of $(v,vr_i)$ using the rotation system defined is $(vr^{o(r_i)-1},v)$, and thus said orbit has $o(r_i)$ elements. In other words, $\left|\{p^m(v,vr_i):m\geq0\}\right|=o(r_i).$ Since $(v,vr_i)\neq (v,vr_j)$ if $i\neq j$, as we vary the values of $i$, the rotation system given by (\ref{eq:rotpmn}) has at least \[
   \sum_{i=1}^n \frac{\left|V(\pmn)\right|}{o(r_i)} = \left|V(\pmn)\right| \sum_{i=1}^n \frac{1}{o(r_i)}
\] regions, where
\[
\sum_{i=1}^n \frac{1}{o(r_i)}=\begin{cases}
    \frac{n}{m},&\text{ if }m\text{ is even,}\\
    \frac{n+1}{2m},&\text{ if }m\text{ is odd.}
\end{cases}
\] 
     
%     \begin{equation}
%         p_v = \begin{cases}
%             (vr_1, vr_2, v\flop_1, v\flop_2, \ldots) & \text{ if } v\in V_1, \text{ and }\\
%             (vr_2, vr_1, v\flop_2, v\flop_1, \ldots), & \text{ if } v\in V_2,
%         \end{cases}
%     \end{equation}
%     where the sets $V_1$ and $V_2$ are assigned according to Algorithm~\ref{alg:labels}.\hl{[DO WE NEED ALL POSSIBLE COMBINATION OF CYCLES OF THE FORM $(r_1^{\pm}r^{\pm}_2)^\ell$. THAT WAS MY IMPRESSION FROM PROP 3.3, BUT IT SEEMS NOT ALL OF THOSE ARE CONSIDERED?]}
%     \hlc[orange]{I believe I have addressed this in Prop 3.3 by appealing to vertex transitivity. Is that sufficient?}

 %   One can verify that the orbit of edges of the form $(v, vr_2)$, if $v \in V_1$ will traverse the $(r_1r_2)^k$ cycle, where $k=2j$ when $m=3j$ or $k=2m$ when $m \neq 3j$.

    Notice that the number of vertices of $\pmn$ is $m^n\cdot n!$ and there are $m^n \cdot n! \cdot n$ edges. Thus 
    \[ 2 - 2 \gamma(\pmn) \geq m^n n! - m^n n! n + \sum_{i=1}^nm^n n!/o(r_i).\] Solving for $\gamma(\pmn)$ yields
    \[
\gamma(\pmn)\leq
\begin{cases}
    \frac{1}{2}m^{n-1} (mn-m-n)n!+1,&\text{ if }m\text{ is even}\\
    \frac{1}{2}m^{n-1} (2 m n - 2 m - n-1) n!+1,&\text{ if }m\text{ is odd,}
\end{cases}
    \] which is the desired result.
       % In the case where $m$ is a multiple of three, there are at least $m^n\cdot n!/(4m/3) = 3/4 m^{n-1} \cdot n!$ distinct orbits of the form $(r_1r_2)^{2m/3}$, and thus there are at least that many regions in the embedding of $\pmn$. Therefore, we have \[ 2 - 2 \gamma(\pmn) \geq m^n \cdot n! - m^n \cdot n! \cdot n + 3/4m^{n-1}\cdot n!.\] Solving for $\gamma(\pmn)$, we see the upper bound of 
    % \[ \gamma(\pmn) \leq \frac{1}{8} m^{n-1}(4mn-4m-3)n! + 1.\]
    % In the case where $m$ is not a multiple of three, there at least $m^n \cdot n!/(4m) = 1/4 m^{n-1} \cdot n!$ distinct orbits of the form $(r_1r_2)^{2m}$, and thus there are at least that many regions in the embedding of $\pmn$. Therefore, we have \[ 2 - 2\gamma(\pmn) \geq m^n \cdot n! - m^n \cdot n! \cdot n + 1/4 m^{n-1} \cdot n!.\] Solving for $\gamma(\pmn)$, we see the upper bound of \[ \gamma(\pmn) \leq \frac{1}{8} m^{n-1}(4mn-4m-1)n! + 1,\] as desired.
\end{proof}

\begin{rem}
    It is of note that the bounds obtained in Theorems~\ref{t:pnupperbound},~\ref{t:bpnupperbound}, and ~\ref{t:pmnupperbound} are asymptotically tight since they have the same asymptotic behavior as the lower bounds given in Remark~\ref{rem:pnbounds} and Corollaries~\ref{cor:bpnlowerbound} and~\ref{cor:pmnlowerbound}. In other words, $\gamma(\bpn)$ is $\Theta(2^nnn!)$ and $\gamma(\pmn)$ is $\Theta(m^nnn!)$.
\end{rem}

We present several values of our upper and lower bounds for $\gamma(\pmn)$ and several values of $m,n$ in Table~\ref{tab:pmn}.

\begin{table}[ht!]
\begin{tabular}{|c|c|r|r|}
\hline
  &  & \multicolumn{2}{|c|}{$\gamma(\pmn)$} \\ \hline
$m$ & $n$ & lower bound & upper bound \\ \hline
3   & 3   & 1 &  217           \\ \hline
4   & 3   & 97 & 241            \\ \hline
5   & 3   & 301 & 1201            \\ \hline
6   & 3   & 649 & 973            \\ \hline
3   & 4   & 325 & 4\;213           \\ \hline
4   & 4   & 3\;073 & 6\;145            \\ \hline
5   & 4   & 10\;501 & 37\;501            \\ \hline
6   & 4   & 25\;921 & 36\;289            \\ \hline
3   & 5   & 9\;721 & 87\;481           \\ \hline
4   & 5   & 92\;161 & 168\;961           \\ \hline
5   & 5   & 375\;001 & 1\;275\;001           \\ \hline
6   & 5   & 1\;088\;641 & 1\;477\;441           \\ \hline
\end{tabular}%
\caption{Values of our bounds for $\gamma(\pmn)$ for several values of $m,n$ with $n>2$. While the bound of Theorem~\ref{t:pmnupperbound} holds for $n=2$, we obtain a better bound in Remark~\ref{r:pm2}.}
\label{tab:pmn}
\end{table}

\begin{rem}\label{r:pm2}
    In the particular case of $n=2$, i.e., $\mathbb{P}_m(2)$ and the rotation system described in the proof of Theorem~\ref{t:pmnupperbound} the authors computed the set of orbits for any $m \geq 3$. With this particular rotation system the boundary cycles are of the form $(r_1)^m$, $(r_2)^a$, and $(\flop_1\flop_2)^{b}$ where 
    \begin{align*}
        a= \begin{cases}
            m, &\text{ if } m \text{ is even}\\
            2m, &\text{ if } m \text{ is odd}
        \end{cases}
        & \qquad b = \begin{cases}
            2m/3, &\text{ if } m = 3j\\
            2m, &\text{ if } m \neq 3j,
        \end{cases}
    \end{align*}
    for $j$ a positive integer. For each vertex $v\in\mathbb{P}_m(2)$ there is only one boundary cycle of the first two forms and two boundary cycles of the last form, one for the $\flop_1$-edge incident to $v$ and the other for the $\flop_2$-edge. The number of vertices in $\mathbb{P}_m(2)$ is $2m^2$ and the number of edges is $4m^2$. We can then count the number of regions/orbits of the rotation system, which there are $2m^2/m + 2m^2/a + 2\cdot2m^2/b$. Combining all of the cases with the Euler-Poincar\'e formula we have the following integer-valued polynomial upper bounds for the genus of $\mathbb{P}_m(2)$,
    \[
        \gamma(\mathbb{P}_m(2)) \leq \begin{cases}
            m^2-3m+1, &\text{ if } m \text{ odd and } m=3j\\
            m^2-2m+1, &\text{ if } m \text{ odd and } m\neq3j\\
            m^2-(7/2) m+1, &\text{ if } m \text{ even and } m=3j\\
            m^2-(5/2) m+1, &\text{ if } m \text{ even and } m\neq3j,\\
        \end{cases}
    \]
    for any positive integer $j$. 

    Specializing even further, when $m=3$, i.e., in the graph $\mathbb{P}_3(2)$ in Figure~\ref{f:p32}, the upper bound is one. Furthermore, one can verify that $\mathbb{P}_3(2)$ contains a $K_{3,3}$ minor (e.g.,\,with $1^02^0, 2^21^1, 1^12^1$ in one independent set and $1^22^0, 2^11^1, 2^21^2$ in the other) and, thus, is not planar. 
\end{rem}

\begin{cor}
    The graph $\mathbb{P}_3(2)$ is toroidal. That is, $\gamma(\mathbb{P}_3(2)) =1.$
\end{cor}

In Table~\ref{tab:pm2}, we compute values of $\gamma(\mathbb{P}_m(2))$ for several values of $m$.

\begin{table}[ht!]
\begin{tabular}{|c|r|r|}
\hline
    & \multicolumn{2}{|c|}{$\gamma(\mathbb{P}_m(2))$} \\ \hline
$m$ & $\lceil$lower bound$\rceil$ & upper bound \\ \hline
3   & -2 &  1            \\ \hline
4   & 1  &  7            \\ \hline
5   & 6  &  16            \\ \hline
6   & 13 &  16           \\ \hline
7   & 18 &  36           \\ \hline
8   & 23 &  45           \\ \hline
9   & 28 &  55           \\ \hline
10  & 35 &  76           \\ \hline
11  & 42 &  100           \\ \hline
12  & 49 &  103           \\ \hline

\end{tabular}%
\caption{Values of our lower bound and specialized upper bound for $\gamma(\mathbb{P}_m(2))$, from Remark~\ref{r:pm2}, for several values of $m$.}
\label{tab:pm2}
\end{table}

%\subsection{Labeling $\pmn$}

%%%%%%%%%%%%%%%%%%%%%%%%%%%%%%%%%%%%%%%%%%%%%%%%%%%%%%%%%%%%%%%%%%%%%%%
% Graph of P_3(2)
\pgfplotsset{compat=1.15}
\usetikzlibrary{arrows}
\definecolor{qqttqq}{rgb}{1.,0.,0.}%red above
\definecolor{qqwuqq}{rgb}{1.,0.,0.}%red
\definecolor{ccqqqq}{rgb}{0.,0.,1.}%blue
\definecolor{ududff}{rgb}{0.,0.,0.}%black
\begin{figure}
    \centering
\begin{tikzpicture}[scale=0.7,line cap=round,line join=round,>=triangle 45,x=1cm,y=1cm, scale=1.2]
\draw [line width=1.5pt,color=ccqqqq] (-12.,6.) -- (-10.,6.);
\draw [line width=1.5pt,color=ccqqqq] (-10.,6.) -- (-11.,8.);
\draw [line width=1.5pt,color=ccqqqq] (-11.,8.) -- (-12.,6.);
\draw [line width=1.5pt,color=qqwuqq] (-14.,7.) -- (-16.,6.);
\draw [line width=1.5pt,color=ccqqqq] (-14.,7.) -- (-13.,9.);
\draw [line width=1.5pt,color=ccqqqq] (-13.,9.) -- (-15.,9.);
\draw [line width=1.5pt,color=ccqqqq] (-15.,9.) -- (-14.,7.);
\draw [line width=1.5pt,color=ccqqqq] (-16.,6.) -- (-17.,8.);
\draw [line width=1.5pt,color=ccqqqq] (-17.,8.) -- (-18.,6.);
\draw [line width=1.5pt,color=ccqqqq] (-18.,6.) -- (-16.,6.);
\draw [line width=1.5pt,color=ccqqqq] (-16.,4.) -- (-18.,4.);
\draw [line width=1.5pt,color=ccqqqq] (-18.,4.) -- (-17.,2.);
\draw [line width=1.5pt,color=ccqqqq] (-17.,2.) -- (-16.,4.);
\draw [line width=1.5pt,color=ccqqqq] (-14.,3.) -- (-15.,1.);
\draw [line width=1.5pt,color=ccqqqq] (-15.,1.) -- (-13.,1.);
\draw [line width=1.5pt,color=ccqqqq] (-13.,1.) -- (-14.,3.);
\draw [line width=1.5pt,color=ccqqqq] (-12.,4.) -- (-11.,2.);
\draw [line width=1.5pt,color=ccqqqq] (-11.,2.) -- (-10.,4.);
\draw [line width=1.5pt,color=ccqqqq] (-10.,4.) -- (-12.,4.);
\draw [line width=1.5pt,color=qqttqq] (-10.,6.) -- (-10.,4.);
\draw [line width=1.5pt,color=qqttqq] (-10.,4.) -- (-17.,8.);
\draw [line width=1.5pt,color=qqttqq] (-17.,8.) -- (-15.,9.);
\draw [line width=1.5pt,color=qqttqq] (-15.,9.) -- (-15.,1.);
\draw [line width=1.5pt,color=qqttqq] (-15.,1.) -- (-17.,2.);
\draw [line width=1.5pt,color=qqttqq] (-17.,2.) -- (-10.,6.);
\draw [line width=1.5pt,color=qqttqq] (-11.,2.) -- (-13.,1.);
\draw [line width=1.5pt,color=qqttqq] (-13.,1.) -- (-13.,9.);
\draw [line width=1.5pt,color=qqttqq] (-13.,9.) -- (-11.,8.);
\draw [line width=1.5pt,color=qqttqq] (-11.,8.) -- (-18.,4.);
\draw [line width=1.5pt,color=qqttqq] (-18.,4.) -- (-18.,6.);
\draw [line width=1.5pt,color=qqttqq] (-18.,6.) -- (-11.,2.);
\draw [line width=1.5pt,color=qqwuqq] (-12.,6.) -- (-14.,7.);
\draw [line width=1.5pt,color=qqwuqq] (-16.,6.) -- (-16.,4.);
\draw [line width=1.5pt,color=qqwuqq] (-16.,4.) -- (-14.,3.);
\draw [line width=1.5pt,color=qqwuqq] (-14.,3.) -- (-12.,4.);
\draw [line width=1.5pt,color=qqwuqq] (-12.,4.) -- (-12.,6.);
\begin{scriptsize}
\draw [fill=ududff] (-12.,6.) circle (2.5pt);
\draw[color=ududff] (-12.470927732826164,5.960200078854992) node {${2^21^0}$};
\draw [fill=ududff] (-10.,6.) circle (2.5pt);
\draw[color=ududff] (-9.615837700803436,6.219133048086954) node {${2^01^0}$};
\draw [fill=ududff] (-11.,8.) circle (2.5pt);
\draw[color=ududff] (-10.779022528664548,8.31678232750867) node {${2^11^0}$};
\draw [fill=ududff] (-13.,1.) circle (2.5pt);
\draw[color=ududff] (-12.727734772743553,0.6600159486987923) node {$\boxed{2^21^1}$};
\draw [fill=ududff] (-17.,8.) circle (2.5pt);
\draw[color=ududff] (-17.33515519479081,8.237632) node {$\circled{2^21^2}$};
\draw [fill=ududff] (-14.,3.) circle (2.5pt);
\draw[color=ududff] (-14.087301454659137,3.427661162009757) node {${2^01^1}$};
\draw [fill=ududff] (-15.,1.) circle (2.5pt);
\draw[color=ududff] (-15.326017764848892,0.6146970593016062) node {$\circled{2^11^1}$};
\draw [fill=ududff] (-16.,6.) circle (2.5pt);
\draw[color=ududff] (-15.524100362355025,5.929987485923535) node {${2^11^2}$};
\draw [fill=ududff] (-18.,6.) circle (2.5pt);
\draw[color=ududff] (-18.468127429720465,6.152963786963938) node {${2^01^2}$};
\draw [fill=ududff] (-17.,2.) circle (2.5pt);
\draw[color=ududff] (-17.35026149125654,1.6872441083683439) node {${1^22^2}$};
\draw [fill=ududff] (-16.,4.) circle (2.5pt);
\draw[color=ududff] (-15.619206658820753,4.313613764090564) node {${1^02^2}$};
\draw [fill=ududff] (-18.,4.) circle (2.5pt);
\draw[color=ududff] (-18.447277057994635,3.966168945378804) node {${1^12^2}$};
\draw [fill=ududff] (-13.,9.) circle (2.5pt);
\draw[color=ududff] (-12.939222923263754,9.37422308010968) node {$\circled{1^22^0}$};
\draw [fill=ududff] (-10.,4.) circle (2.5pt);
\draw[color=ududff] (-9.630943997269165,3.66404301606423) node {$\boxed{1^12^1}$};
\draw [fill=ududff] (-15.,9.) circle (2.5pt);
\draw[color=ududff] (-14.99367924260286,9.37422308010968) node {$\boxed{1^02^0}$};
\draw [fill=ududff] (-11.,2.) circle (2.5pt);
\draw[color=ududff] (-10.733703639267361,1.807669294231259) node {${1^02^1}$};
\draw [fill=ududff] (-14.,7.) circle (2.5pt);
\draw[color=ududff] (-14.102407751124865,6.6910463844700715) node {${1^12^0}$};
\draw [fill=ududff] (-12.,4.) circle (2.5pt);
\draw[color=ududff] (-12.470927732826164,4.35893265348775) node {${1^22^1}$};
\end{scriptsize}
\end{tikzpicture}
    \caption{The graph $\mathbb{P}_3(2)$ contains $K_{3,3}$ as a minor where the  vertices forming the two independent sets are inside a rectangle or a circled rectangle, respectively. Therefore $\mathbb{P}_3(2)$ is not planar.}
    \label{f:p32}
\end{figure}
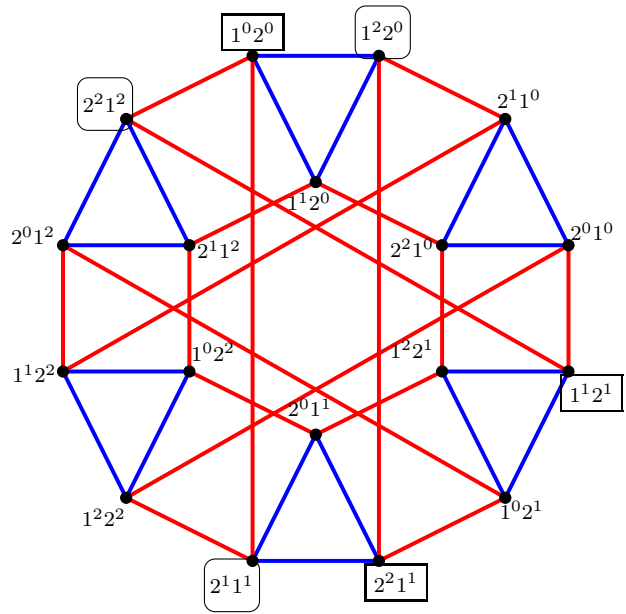

\section{Concluding remarks}

Our main contribution is determining upper bounds for the genus of $\pn,\bpn$ and $\pmn$ and lower bounds for the genus of $\bpn$ and $\pmn$. Our upper bound for $\gamma(\pn)$ is the tightest bound known in the literature, and the bounds for $\gamma(\bpn)$ and $\gamma(\pmn)$ are the first-known bounds. A key ingredient in our proofs is finding rotation systems where certain cycles in the graphs considered become boundaries for regions of a 2-cell embedding. In the case of $\gamma(\pmn)$, we took advantage of cycles given by subsequently traversing the same prefix reversal. However, in the case of $\pn$ and $\bpn$, since the prefix reversals are involutions (their order is 2), we needed to employ a different technique, using cycles with two distinct reversals. Indeed, we devised a labeling algorithm, called $\AlGra$, that we used to established that certain cycles in the pancake graph $\pn$ and burnt pancake graph $\bpn$ become boundaries of regions of a 2-cell embedding given by a particular rotation system. We remark that there are more orbits in the rotation systems used and that in principle one could sharpen the bounds found. However, we were unable to determine a general patter for \emph{all} the orbits in our rotation systems. It is also possible that there are rotation systems different from the ones we used that would produce a lower genus. However, after a computer aided application of different rotation systems for $\mathbb{P}_4$, we obtained the same number of orbits/regions as the number that we obtained using the rotation system employed in our proofs. So there is some evidence to suggest that our rotation systems may be optimal. 

\section*{Acknowledgment} The authors wish to thank Richard Ehrenborg for his motivating question regarding the genus of $\mathbb{P}_4$ during a talk the first author gave at the Combinatorics seminar at the University of Kentucky.

%\bibliography{bib}

\begin{thebibliography}{10}

\bibitem{BjornerBrenti}
Anders Bj{\"o}rner and Francesco Brenti.
\newblock {\em Combinatorics of {C}oxeter groups}, volume 231 of {\em Graduate
  Texts in Mathematics}.
\newblock Springer, New York, 2005.

\bibitem{BBP19}
Sa\'{u}l~A. Blanco, Charles Buehrle, and Akshay Patidar.
\newblock Cycles in the burnt pancake graph.
\newblock {\em Discrete Appl. Math.}, 271:1--14, 2019.

\bibitem{BBrelations}
Sa\'{u}l~A. Blanco and Charles Buehrle.
\newblock Some relations on prefix-reversal generators of the symmetric and
  hyperoctahedral groups.
\newblock {\em Australas. J. Combin.}, 76(part 1):404--427, 2020.

\bibitem{BB22}
Sa\'ul~A. Blanco and Charles Buehrle.
\newblock Presentations of {C}oxeter groups of type {$A$}, {$B$}, and {$D$}
  using prefix-reversal generators.
\newblock {\em Applicable Algebra in Engineering, Communication and Computing},
  2022.

\bibitem{BB23}
Sa\'{u}l~A. Blanco and Charles Buehrle.
\newblock Lengths of cycles in generalized pancake graphs.
\newblock {\em Discrete Mathematics}, 346(12):113624, 2023.




\bibitem{CSW21}
Ben Cameron, Joe Sawada, and Aaron Williams.
\newblock A {H}amilton {C}ycle in the $k$-sided {P}ancake {N}etwork.
\newblock In Paola Flocchini and Lucia Moura, editors, {\em Combinatorial
  Algorithms - 32nd International Workshop, {IWOCA} 2021, Ottawa, ON, Canada,
  July 5-7, 2021, Proceedings}, volume 12757 of {\em Lecture Notes in Computer
  Science}, pages 137--151. Springer, 2021.

\bibitem{CLZ16}
Gary Chartrand, Linda Lesniak, and Ping Zhang.
\newblock {\em Graphs \& {D}igraphs}.
\newblock Textbooks in Mathematics. CRC Press, Boca Raton, FL, sixth edition,
  2016.

\bibitem{C07}
Jianer Chen, Iyad~A. Kanj, Ljubomir Perkoviƒá, Eric Sedgwick, and Ge~Xia.
\newblock Genus characterizes the complexity of certain graph problems: Some
  tight results.
\newblock {\em Journal of Computer and System Sciences}, 73(6):892--907, 2007.

\bibitem{Dweighter77}
Harry Dweighter, Michael~R. Garey, David~S. Johnson, and Shen Lin.
\newblock Problems and {S}olutions: {S}olutions of {E}lementary {P}roblems:
  {E}2569.
\newblock {\em Amer. Math. Monthly}, 84(4):296, 1977.

\bibitem{D88}
Walther Dyck.
\newblock Beitr\"{a}ge zur {A}nalysis situs.
\newblock {\em Math. Ann.}, 32(4):457--512, 1888.

\bibitem{E60}
John Robert~Jr. Edmonds.
\newblock A {C}ombinatorial {R}epresentation for {O}riented {P}olyhedral
  {S}urfaces.
\newblock Master's thesis, University of Maryland, College Park, 1960.

\bibitem{EL22}
Louis Esperet and Benjamin L\'ev\^eque.
\newblock Local certification of graphs on surfaces.
\newblock {\em Theoretical Computer Science}, 909, 01 2022.

\bibitem{Fetal20}
Laurent Feuilloley, Pierre Fraigniaud, Pedro Montealegre, Ivan Rapaport,
  \'{E}ric R\'{e}mila, and Ioan Todinca.
\newblock Compact distributed certification of planar graphs.
\newblock In {\em Proceedings of the 39th Symposium on Principles of
  Distributed Computing}, PODC '20, page 319--328, New York, NY, USA, 2020.
  Association for Computing Machinery.

\bibitem{Fetal23}
Laurent Feuilloley, Pierre Fraigniaud, Pedro Montealegre, Ivan Rapaport, Eric
  Remila, and Ioan Todinca.
\newblock Local certification of graphs with bounded genus.
\newblock {\em Discrete Applied Mathematics}, 325:9--36, 01 2023.

\bibitem{GatesPapa}
William~H. Gates and Christos~H. Papadimitriou.
\newblock Bounds for sorting by prefix reversal.
\newblock {\em Discrete Math.}, 27(1):47--57, 1979.

\bibitem{HannenPev}
Sridhar Hannenhalli and Pavel~A. Pevzner.
\newblock Transforming cabbage into turnip: polynomial algorithm for sorting
  signed permutations by reversals.
\newblock {\em J. ACM}, 46(1):1--27, 1999.

\bibitem{Haynes2008}
Karmella~A. Haynes, Marian~L. Broderick, Adam~D. Brown, Trevor~L. Butner,
  James~O. Dickson, W.~Lance Harden, Lane~H. Heard, Eric~L. Jessen, Kelly~J.
  Malloy, Brad~J. Ogden, Sabriya Rosemond, Samantha Simpson, Erin Zwack,
  A.~Malcolm Campbell, Todd~T. Eckdahl, Laurie~J. Heyer, and Jeffrey~L. Poet.
\newblock Engineering bacteria to solve the {B}urnt {P}ancake {P}roblem.
\newblock {\em Journal of Biological Engineering}, 2(1):8, May 2008.

\bibitem{H91}
L.~Heffter.
\newblock Ueber das {P}roblem der {N}achbargebiete.
\newblock {\em Math. Ann.}, 38(4):477--508, 1891.

\bibitem{H98}
L.~Heffter.
\newblock Ueber metacyklische {G}ruppen und {N}achbarconfigurationen.
\newblock {\em Math. Ann.}, 50(2-3):261--268, 1898.

\bibitem{KF95}
Arkady Kanevsky and Chao Feng.
\newblock On the embedding of cycles in pancake graphs.
\newblock {\em Parallel Computing}, 21(6):923 -- 936, 1995.

\bibitem{Dweighter75}
D.~J. Kleitman, Edvard Kramer, J.~H. Conway, Stroughton Bell, and Harry
  Dweighter.
\newblock Problems and {S}olutions: {E}lementary {P}roblems: {E}2564-{E}2569.
\newblock {\em Amer. Math. Monthly}, 82(10):1009--1010, 1975.

\bibitem{KKP05}
Amos Korman, Shay Kutten, and David Peleg.
\newblock Proof labeling schemes.
\newblock {\em Distributed Computing}, 22:215--233, 07 2005.

\bibitem{K30}
Casimir Kuratowski.
\newblock Sur le probl\`eme des courbes gauches en topologie.
\newblock {\em Fundamenta Mathematicae}, 15(1):271--283, 1930.

\bibitem{LJD93}
S.~Lakshmivarahan, Jung-Sing Jwo, and S.K. Dhall.
\newblock Symmetry in interconnection networks based on {C}ayley graphs of
  permutation groups: A survey.
\newblock {\em Parallel Computing}, 19(4):361 -- 407, 1993.

\bibitem{NB09}
Quan~T. Nguyen and Said Bettayeb.
\newblock The upper bound and lower bound of the genus of pancake graphs.
\newblock In {\em 2009 IEEE Symposium on Computers and Communications}, pages
  800--804, 2009.

\bibitem{R65}
Gerhard Ringel.
\newblock Das {G}eschlecht des vollst\"{a}ndigen paaren {G}raphen.
\newblock {\em Abh. Math. Sem. Univ. Hamburg}, 28:139--150, 1965.

\bibitem{T89}
Carsten Thomassen.
\newblock The graph genus problem is {NP}-complete.
\newblock {\em Journal of Algorithms}, 10(4):568--576, 1989.

\bibitem{W37}
K.~Wagner.
\newblock \"Uber eine eigenschaft der ebenen komplexe.
\newblock {\em Mathematische Annalen}, 114:570--590, 1937.

\bibitem{W22}
Liangxia Wan.
\newblock The genus of a graph: A survey.
\newblock {\em Symmetry}, 15(2), 2023.

\bibitem{West00}
Douglas~B. West.
\newblock {\em Introduction to Graph Theory}.
\newblock Prentice Hall, 2 edition, September 2000.

\bibitem{W72}
Arthur~T. White.
\newblock On the genus of a group.
\newblock {\em Trans. Amer. Math. Soc.}, 173:203--214, 1972.

\bibitem{Y63}
J.~W.~T. Youngs.
\newblock Minimal imbeddings and the genus of a graph.
\newblock {\em J. Math. Mech.}, 12:303--315, 1963.

\end{thebibliography}

\addresseshere

\end{document}